\documentclass[a4paper,11pt]{amsart}

\usepackage{amssymb,amsmath,stmaryrd,mathrsfs, amscd, extarrows, mathtools}

\usepackage{graphics,color,graphicx}

\definecolor{darkblue}{RGB}{0,0,160}
\definecolor{darkred}{RGB}{180,0,0}

\usepackage[extension=pdf,pdfencoding=auto,unicode=true,bookmarks=true,linktocpage=true,pdftex]{hyperref}
\usepackage{bookmark}
\hypersetup{colorlinks, 
allcolors=darkblue,
pdftex
}

\hypersetup{bookmarks=true}
\hypersetup{
pdfauthor={Fabio Trova}
}
\hypersetup{
pdftitle={Nakayama categories and groupoid quantization}
}

\input xy
\xyoption{all}
\xyoption{2cell}
\UseAllTwocells
\CompileMatrices
\UseTwocells
\xyoption{arc}

\DeclareRobustCommand{\SkipTocEntry}[5]{}

\usepackage{cleveref}

\newtheorem{thm}{Theorem}[section]
\crefname{thm}{Theorem}{Theorems}
\newtheorem{cor}[thm]{Corollary}
\newtheorem{prop}[thm]{Proposition}
\crefname{prop}{Proposition}{Propositions}
\newtheorem{lem}[thm]{Lemma}
\crefname{lem}{Lemma}{Lemmas}

\newtheorem{defin}[thm]{Definition}
\crefname{defin}{Definition}{Definitions}
\theoremstyle{definition}
\newtheorem{rem}[thm]{Remark}
\crefname{rem}{Remark}{Remarks}
\newtheorem{ex}[thm]{Example}
\crefname{ex}{Example}{Examples}

\newcommand{\cat}[1]{\mathcal{#1}} 

\newcommand{\defeq}{\mathrel{\mathop:}=}  

\newcommand{\sam}{\ensuremath{\mathit{Sum}}}  
\newcommand{\pro}{\ensuremath{\mathit{Prod}}}  

\newcommand{\map}[2]{\mathit{Map}(#1,#2)}   

\newcommand{\Hom}[3]{#1(#2,#3)}   

\newcommand{\lad}{\dashv}  

\newcommand{\iso}{\cong}  

\newcommand{\eps}{\varepsilon} 

\newcommand{\vect}{\mathit{Vect}} 

\DeclareMathOperator*{\tensor}{\otimes}  
\DeclareMathOperator*{\etensor}{\boxtimes}  

\newcommand{\ra}[1]{{#1}_{*}}   
\newcommand{\la}[1]{{#1}_{!}}   

\newcommand{\rma}[1]{{#1}_{R}}  
\newcommand{\lma}[1]{{#1}_{L}}  

\newcommand{\etar}[1]{\ra{\eta}^{#1}} 
\newcommand{\epsr}[1]{\ra{\eps}^{#1}} 

\newcommand{\etal}[1]{\la{\eta}^{#1}} 
\newcommand{\epsl}[1]{\la{\eps}^{#1}} 

\newcommand{\ga}[1]{\gamma_{#1}} 
\newcommand{\hga}[1]{\hat{\gamma}_{#1}} 
\newcommand{\na}[1]{\nu_{#1}} 

\newcommand{\xto}[1]{\xrightarrow{#1}}

\newcommand{\famic}[2]{\ensuremath{\mathrm{Fam}_{#1}(#2)}}  

\newcommand{\fami}[1]{\ensuremath{\mathrm{Fam}_{1}^{1}(#1)}}  

\newcommand{\bordnk}[1]{\mathrm{Bord}_{n}^{#1}} 

\newcommand{\nalg}[1]{\ensuremath{\mathrm{Alg}_{#1}}} 

\newcommand{\fgrp}{\mathrm{Gpd}_{fin}} 

\begin{document}
\title{Nakayama Categories and Groupoid Quantization}
\author{Fabio Trova}
\address{Max-Planck-Institut f\"ur Mathematik, Bonn, Germany}
\email{fabiotrova@mpim-bonn.mpg.de\\ fabio.trova@gmail.com}
\date{\today}

\begin{abstract}
We  provide a precise description, albeit in the  situation of standard categories, of the quantization functor $\sam$ proposed by D.S. Freed, M.J. Hopkins, J. Lurie, and C. Teleman in a way enough abstract and flexible to suggest  that an extension of the construction to the general context of ($\infty,n$)-categories  should indeed be possible. Our method is in fact based primarily on dualizability and adjunction data, and is well suited for the homotopical setting. The construction also sheds light on the need  of certain rescaling automorphisms, and in particular on the nature and properties of the Nakayama isomorphism.  
\end{abstract}
\maketitle

\tableofcontents

\section{Introduction}

A $k$-extended $n$-dimensional $\cat{C}$-valued topological quantum field theory (TQFT for short) is a symmetric monoidal functor $Z\colon\bordnk{n-k}\to\cat{C}$, from the $(\infty,k)$-category of $n$-cobordisms extended down to codimension $k$ \cite{cs,lu3} to a symmetric monoidal $(\infty,k)$-category $\cat{C}$. For $k=1$ one recovers the classical TQFT's as introduced in \cite{at}, while $k=n$ gives rise to \emph{fully} extended TQFT's which, by  the Cobordism Hypothesis \cite{bd,lu3}, we know being in correspondence with the fully dualizable objects of $\cat{C}$.

In \cite{fhlt} the authors suggested the existence of an extended TQFT of Dijkgraaf-Witten type, associated to any given representation  ${V\colon X\to\cat{C}}$ of a finite $\infty$-groupoid $X$ with values in a symmetric monoidal $(\infty,k)$-category $\cat{C}$; under the identification of $\infty$-groupoids with homotopy types, $V$ is equivalently seen as the datum of a local system on a topological space $X$ of finite homotopy type, with values in $\mathcal{C}$.  
Such a theory appears as the composition of two monoidal functors
\begin{align*}
\mathit{Bun}_{k,V} &\colon \bordnk{n-k} \to \famic{k}{\cat{C}}  \\
\sam_{k} &\colon \famic{k}{\cat{C}} \to \cat{C}
\end{align*}
where $\famic{k}{\cat{C}}$ is an $(\infty,k)$-category of representations, with objects spaces equipped with a functor to $\cat{C}$ and morphisms the correspondences between these (see \cite{ha2,lu3}). 
In the fully extended case the first functor associates to the point the given representation $V$, while $\sam$ is the proper quantization functor carrying a representation $W\colon Y\to\cat{C}$ to its colimit. This can be thought of as taking the global sections of the local system $W$, or, from a more algebraic point of view, as the space of coinvariants for the action of the $\infty$-groupoid $Y$ on $W$. 

Unfortunately in \cite{fhlt} not many details are provided on how the above  machinery should behave; the article in fact only addresses the above colimit construction and demands that colimits and limits in $\cat{C}$ coincide in order to get a well defined functor. A better picture is given in the case when $\cat{C}$ is the $(\infty,n)$-category $\nalg{n-1}$ (in the sense of \cite[Def. 7.1]{fhlt}) of iterated algebras and modules over $\vect$, and $V \colon X \to \nalg{n-1}$ is the trivial representation collapsing everything on the unit.  
In this situation, a detailed definition of the functor $\sam$ has been given for $n=2$ in \cite{mo1}, while in \cite{mo2} also spaces endowed with a $2$-cocycle have been considered; again, the author requires the existence of an isomorphism between right and left adjoints, given explicitly by the Nakayama isomorphism (cfr. \cite{ben}).  
It is worth noting that, despite the Cobordism Hypothesis, neither \cite{fhlt} or \cite{mo1,mo2}  seem to make  use of dualizability in their construction.\\

The primary goal of this paper is to provide a detailed construction of the functor $\sam_{1}$ and therefore of $n$-dimensional theories in the sense of Atiyah. 
The method we propose differs from the ones above in two main points:
\begin{itemize}
\item[-] It does not  assume the existence of any a priori given isomorphism between right and left adjoints. 
\item[-] It relies almost only on existence of duals, in line with the Cobordism Hypothesis.
\end{itemize}

The main obstacles to the definition of the functor $\sam$ appear when studying its behaviour at the level of morphisms. 
In fact, to a morphism $\alpha$ in $\famic{n}{\cat{C}}$ we could associate a morphism $\sam(\alpha)$ in $\cat{C}$ via a pull-push procedure, as soon as we have a natural transformation (which is \emph{not} required to be invertible) $\ra{f}\to\la{f}$ from the right to the left adjoint of the pullback functor $f^*$, for any map $f\colon X\to Y$ of spaces.

A candidate $\ga{f}$ is presented in \cref{ga}; the map was originally defined in \cite{ha1} as long as certain projection formulas were satisfied, which actually always happens (\cref{lproj,rproj}) in presence of duals.  
Despite being very well behaved with respect to homotopy squares and external tensor products (as shown in \cref{prenakasubsec,prenakakantensorsubsec}), the map $\ga{}$ still lacks an essential feature. In fact, in order to make $\sam$ into a functor, we need the assignment $f\mapsto\ga{f}$ to be functorial as well, in the sense that for any two composable maps $X\xto{f} Y \xto{g} Z$ the following triangle
\[
\xymatrix@!=5pt{
& \la{g}\ra{f} \ar[rd]^{\la{g}\ga{f}} & \\
\ra{g}\ra{f} 
\ar[ur]^{\ga{g}} \ar[rr]_{\ga{gf}} && \la{g}\la{f} \\
}
\]
should commute.

It can be seen that such condition is generally not verified, due to the possible non-triviality of the fibers of the maps $f,g$ involved. As shown in \cref{counterexcn,counterexcncm}, when the target category $\cat{C}$ is that of finite vector spaces over a positive characteristic field, one of the two paths might be the zero map while the other one is nonzero. Instead,  in characteristic zero, the two ways of going from the right to the left adjoints differ by an invertible scalar. 

One might then hope  the assignment $f\mapsto\ga{f}$ to define a ``projective functor'', and the discrepancy to result in a ``cocycle condition''. 
This is actually the situation in the case of characteristic zero fields: we can in fact rescale the above map $\ga{}$ by a factor $\delta_{f}^{-1}$ related to the Baez-Dolan groupoid cardinality \cite{bd2}, so to neutralize the anomaly.  
In categorical terms, $\delta_{f}$ can be described via $\ga{}$ itself (\cref{weightdelta}) and we can define a new map (the \emph{Nakayama map}, \cref{naka})
\[
\na{f}\colon \ra{f}\to\la{f}
\]
by multiplying $\ga{f}$ by $\delta_{f}^{-1}$. 
The best one can now achieve is to ask this new map to be functorial in the above sense.  
We will therefore define a \emph{Nakayama category} (\cref{nakacontext}) as a rigid category, in which for any map $f$ of spaces the endomorphism $\delta_{f}$ is invertible and the assignment $f\mapsto\na{f}$ is functorial.  

Under these hypotheses in \cref{thmsum,thmprod} we will construct, by means of Kan extensions, two monoidal quantization functors 
\begin{align*}
\sam &\colon \famic{1}{\cat{C}} \to \cat{C} \\
\pro &\colon \famic{1}{\cat{C}} \to \cat{C}
\end{align*}
The map $\na{}$ will then define a monoidal transformation between them and turn out to be an isomorphism (\cref{nakamontran,nakaiso}), thus answering an ambidexterity problem \cite{hel,hol} in the context of categories with duals. 

As one would expect, when $\cat{C}$  is the category of finite dimensional vector spaces over a characteristic zero field, the map $\na{f}$ recovers the known Nakayama isomorphism, and  justifies the choices made in \cite{fhlt,mo1,mo2} as instances of a more general and abstract framework. In particular, it explains from a new point of view the obstruction to quantizing over fields of positive characteristic. \\

Although the present paper only deals with standard categories and groupoids, the methods and ideas we propose are intrinsically homotopical, and formulated mainly in terms of adjunction and duality data via calculus of mates and Beck-Chevalley conditions. This suggests it should be possible to export them to the framework of higher categories so to build the proper $n$-dimensional theories which are the true aim of \cite{fhlt}, a topic which will be investigated elsewhere.

\subsection{Conventions}
While most of the constructions can be stated in the general case of non-symmetric monoidal categories, we will be interested into symmetric ones. By \emph{monoidal category} we will therefore mean \emph{symmetric} monoidal category.  
By space we will always mean $\infty$-groupoid; in particular, also groupoids (i.e. homotopy $1$-types) will be often referred to as spaces.


\addtocontents{toc}{\SkipTocEntry}
\subsection*{Aknowledgements}
I am happily in debt with Domenico Fiorenza, who supervised my PhD thesis from which part of this work is extracted; I'm glad to thank him for the many and helpful discussions, during and after my PhD.  I would like to thank as well Claudia Scheimbauer, Urs Schreiber and Christoph Schweigert for  interesting and fruitful conversations and comments.

\section{Preliminaries} \label{secpre}

We quickly review here below some notions and results, which are essential to the development of our machinery. 

\subsection{Duals}

Recall first of all the notion of dualizable object in a monoidal category $\cat{C}$, which generalises and axiomatises the familiar concept of dual of a vector space. Dualizability is a key tool in this work, and more generally the main and possibly only ingredient to the construction of an extended Topological Quantum Field Theory, which is the ultimate goal of this paper.
\begin{defin} \label{dualcat}
Let $\cat{C}$ be a monoidal category and $A$ an object of $\cat{C}$. 
An object $A^{d}$ of $C$ is a \emph{dual} of $A$ if there are coevaluation and evaluation maps
\begin{align*}
co: 1\to A^{d}\tensor A \\
ev: A\tensor A^{d}\to 1
\end{align*}
such that the compositions
\begin{align*}
A\iso A\tensor 1 &\xrightarrow{id\tensor co} A\tensor A^{d}\tensor A \xrightarrow{ev\tensor id} 1\tensor A\iso A \\
A^{d} \iso 1\tensor A^{d} &\xrightarrow{co\tensor id} A^{d}\tensor A\tensor A^{d} \xrightarrow{id\tensor ev} 1\tensor A^{d}\iso A^{d}
\end{align*}
both equal the identity. 

We say that an object $A$ is \emph{dualizable} if it has a dual. 
We say that a  monoidal category $\cat{C}$ \emph{has duals} (or is a category with duals), if every object of $\cat{C}$ is dualizable.
\end{defin}
Notice that, being $\cat{C}$ symmetric monoidal, $A$ defines a dual to $A^{d}$.\\

As for the case of $\vect$-valued group representations, not only the monoidal structure is transferred to presheaf categories but also dualizability: 

\begin{prop} \label{fundualcat}  \label{funmoncat}
Let $\cat{C}$ be a monoidal category with duals. Then, for any groupoid $X$, the functor category $\cat{C}^{X}$ is a monoidal category with duals.
\end{prop}

The notion of dual can be nicely reformulated in terms of adjoints, in the following sense
\begin{prop} \label{dualadjoint}
Let $\cat{C}$ be a monoidal category, and $A$ an object of $\cat{C}$.  
If $A$ has a dual $A^{d}$, then the coevaluation and evaluation maps define adjunctions
\begin{align*}
A\tensor - &\lad A^{d}\tensor - \\
-\tensor A^{d} &\lad -\tensor A
\end{align*}
By universal property of the adjoints, it follows that any two duals to the same object $A$ are isomorphic; in particular, one has $(A^{d})^{d}\iso A$.
\end{prop}
The latter is a very convenient way to deal with duals; by treating them as adjoint functors, in fact, we will be able to mix  adjunction and dualizability data and place them in the framework of Beck-Chevalley transformations.

As a toy version of the cobordism hypothesis, finally, one can easily prove the following two results:

\begin{prop} \label{monfctdual}
Let $\cat{C},\cat{D}$ be monoidal categories and $A\in\cat{C}$ an object having a  dual $A^{d}$. Then
\begin{itemize}
\item[i)] For any monoidal functor $F\colon\cat{C}\to\cat{D}$, $FA^{d}$ is a dual to $FA$.
\item[ii)] If $\phi\colon F\to G$ is a monoidal transformation between monoidal functors $F,G\colon\cat{C}\to\cat{D}$, then $\phi_{A}\colon FA\to GA$ is an isomorphism.
\end{itemize}
\end{prop}

\begin{cor} \label{montraiso}
Let $\cat{C},\cat{D}$ be monoidal categories, and assume $\cat{C}$ has duals. Then, any monoidal transformation $\phi\colon F\to G$ of monoidal functors from $\cat{C}$ to $\cat{D}$ is an isomorphism.
\end{cor}

\subsection{Mates and Beck-Chevalley condition}

We recollect now some results on the formation and properties of Beck-Chevalley transformed two-cells. 
This is a very helpful formalism, suited for the homotopical setting as well (see \cite[4.7.5.13 et seqq.]{lu4}), to handle adjunction data and functors, and to deal with Kan extensions.

\begin{defin} 
\label{lmate}\label{lbc}  Let
\[
\xymatrix{
\cat{D} \ar[r]^{K} \ar[d]_{G}  & \cat{D}' \ar[d]^{G'}  \\
\cat{C} \ar[r]_{H}     & \cat{C}' \ultwocell<\omit>{\psi}
}
\]
be a lax square of functors. Assume that $G$ and $G'$ have left adjoints $F,F'$. The \emph{left mate} of $\psi$ is the map $\lma{\psi}$ defined as
\[
F'H\xto{F'H\eta} F'HGF \xto{F'\psi}F'G'KF \xto{\eps'} KF
\]

We say that the diagram satisfies the \emph{left Beck-Chevalley condition} if the left mate of $\psi$ is an isomorphism.
\end{defin}

\begin{defin}
\label{rmate}  \label{rbc} Let
\[
\xymatrix{
\cat{C} \ar[r]^{H} \ar[d]_{F} \drtwocell<\omit>{\phi}  & \cat{C}' \ar[d]^{F'}  \\
\cat{D} \ar[r]_{K}     & \cat{D}' 
}
\]
be an oplax square of functors. Assume that $F$ and $F'$ have right adjoints $G,G'$. The \emph{right mate} of $\phi$ is the map $\phi_{R}$  defined as
\[
HG\xto{\eta'} G'F'HG \xto{G'\phi}G'KFG \xto{G'K\eps} G'K
\]
We say that the diagram satisfies the \emph{right Beck-Chevalley condition} if the right mate of $\phi$ is an isomorphism.
\end{defin}

Particularly useful will be the nice behaviour of the above formalism with respect to both horizontal and vertical pasting of diagrams; the following  lemma, and its dual version, will in fact be (ab)used all through this paper.

\begin{lem}  \label{pasteleft} 
\begin{itemize}
\item[i)]
Let 
\[
\xymatrix{
\cat{D} \ar[r]^{K} \ar[d]_{G}  & \cat{D}' \ar[d]_{G'}   \ar[r]^{K'}   & \cat{D}'' \ar[d]^{G''}    \\
\cat{C} \ar[r]_{H}     & \cat{C}' \ultwocell<\omit>{\psi'}   \ar[r]_{H'}     & \cat{C}'' \ultwocell<\omit>{\psi''}
}
\]
be lax squares, with $G,G',G''$ right adjoints. Let $\psi\colon H'HG\to G''K'K$ be the composition of $\psi'$ and $\psi''$. Then the left mate $\lma{\psi}$ of $\psi$ is the composition of $\psi'_{L}$ and $\psi''_{L}$.
\item[ii)]
Let
\[
\xymatrix{
\cat{E} \ar[r]^{L} \ar[d]_{G_{2}}      & \cat{E}' \ar[d]^{G_{2}'} \\
\cat{D} \ar[r]^{K} \ar[d]_{G_{1}}  & \cat{D}' \ar[d]^{G_{1}'}  \ultwocell<\omit>{\psi''}  \\
\cat{C} \ar[r]_{H}     & \cat{C}' \ultwocell<\omit>{\psi'}
}
\]
be lax squares, with $G_{1},G_{2},G'_{1},G'_{2}$ right adjoints. Let $\psi\colon HG_{1}G_{2}\to G_{1}'G_{2}'L$ be the composition of $\psi''$ and $\psi'$. Then $\lma{\psi}$ is the composition of $\lma{\psi''}$ and $\lma{\psi'}$.
\end{itemize}
\end{lem}

Clearly  a specular statement holds for oplax squares and right mates; we leave to the reader the simple task to write it out.\\

It can be seen, as well, that the left and right Beck-Chevalley conditions interact in some cases and are essentially equivalent; more precisely we get

\begin{prop} \label{leftrightbc}
Let
\[
\xymatrix{
\cat{A} \ar[r]^{F_{3}} \ar[d]_{G_{1}}  & \cat{B} \ar[d]^{G_{2}}  \\
\cat{C} \ar[r]_{F_{4}}     & \cat{D} \ultwocell<\omit>{\theta}
}
\]
be a lax square of functors, with $F_{1}\lad G_{1},F_{2}\lad G_{2},F_{3}\lad G_{3},F_{4}\lad G_{4}$. Then the following are equivalent
\begin{itemize}
\item[i)] The diagram satisfies the left Beck-Chevalley condition
\item[ii)] The diagram satisfies the right Beck-Chevalley condition
\end{itemize}
\end{prop}

In the next pages we will be especially interested in diagrams arising from restriction functors. 
Recall that given $f\colon\cat{A}\to\cat{B}$ and a category $\cat{C}$, one has a restriction functor $f^{*}\colon\cat{C}^{\cat{B}}\to\cat{C}^{\cat{A}}$. Its left and right adjoints, if they exist, are given by Kan extensions  and will be denoted respectively by $\la{f}$ and $\ra{f}$; a sufficient condition to  their existence is given by cocompleteness, resp. completeness, of $\cat{C}$ (see \cite{ml}).

\begin{defin}
Let 
\[
\xymatrix{
\cat{A} \ar[r]^{u} \ar[d]_{p} \drtwocell<\omit>{\phi}  & \cat{A}' \ar[d]^{q}  \\
\cat{B} \ar[r]_{v}     & \cat{B}' 
}
\]
be a (op)lax square of functors. 
We say that the  square is \emph{left} (resp. \emph{right}) \emph{exact} if the corresponding diagram
\[
\xymatrix{
\cat{C}^{\cat{B}'}\ar[r]^{v^{*}} \ar[d]_{q^{*}}  & \cat{C}^{\cat{B}} \ar[d]^{p^{*}}  \\
\cat{C}^{\cat{A}'} \ar[r]_{u^{*}}     & \cat{C}^{\cat{A}} \ultwocell<\omit>{\phi^{*}}
}
\]
satisfies the left (right) Beck-Chevalley condition, for any cocomplete (resp. complete) category $\cat{C}$. 
We say that it is an \emph{exact square} if it is both left and right exact.
\end{defin}

A natural question is whether one can find conditions for a square as above to be exact; a positive answer is given by pullbacks along Grothendieck fibrations, which specializes to the following, pivotal

\begin{prop}\label{grpBC}
Let 
\[
\xymatrix{
P \ar[r] \ar[d] & H \ar[d]^{f} \\
K \ar[r]_{g}  & G \ultwocell<\omit>{\phi}
}
\]
be a homotopy pullback of groupoids. 
Then both the diagrams relative to $\phi$ and $\phi^{-1}$ are exact.
\end{prop}

\section{A comparison map between adjoints}  \label{prenakasec}

The present section focuses on the key ingredient of this paper, that is  a natural morphism relating right and left Kan extensions along maps of groupoids. This is a map $\ga{}$, first appeared in \cite{ha1}, which can be defined for any monoidal functor $f^{*}$ having both left and right adjoints and satisfying certain projection formulas. We'll see that dualizability is a sufficient condition to the existence of $\ga{}$, and via the calculus of mates we'll deduce  several interesting properties exhibiting $\ga{}$ as well behaved w.r.t. homotopy exact squares and external tensor products. The results and constructions  given here will provide the background for the map properly used in our quantization functors, which we'll define and study in \cref{secnaka}. 

The situation we are describing can be placed in the framework of Grothendieck's six operations (see for example \cite{cd}) and, specifically, be understood as a Wirthm\"uller context. In view of the latter a map (\cref{mayga}) similar to $\ga{}$ was previously defined in \cite{fhm}, and in fact the two can be shown to coincide (\cref{gamay}).

\subsection{Projection formulas}

The aforementioned map $\ga{}$ will arise thanks to certain ubiquitous  projection morphisms $\lambda$ and $\rho$. The key condition to existence of $\ga{}$ is invertibility of such maps, which we show here below being always satisfied thanks to dualizability. In this subsection we will moreover establish some useful lemmas concerning the maps $\lambda$ and $\rho$, by means of Beck-Chevalley conditions and transformations.

\begin{prop}[Left Projection Formula]\label{lproj} Let $f^{*}\colon \cat{C}\to\cat{D}$ be a monoidal functor, having a left adjoint $\la{f}$. Then, there is a natural map  
\[
\lambda\colon \la{f}(A\tensor f^{*}B) \to \la{f}A\tensor B
\]
If $\cat{C}$ has duals, then $\lambda$ is an isomorphism. 
\proof
Consider the square
\[
\xymatrix{
\cat{C} \ar[r]^{-\tensor B} \ar[d]_{f^{*}}  & \cat{C} \ar[d]^{f^{*}}  \\
\cat{D} \ar[r]_{-\tensor f^{*}B}     & \cat{D} \ultwocell<\omit>{\mu}
}
\]
and let $\lambda$ be the left mate of $\mu$. 
Invertibility of $\lambda$ is just the left Beck-Chevalley condition for the above diagram.   
Since $\cat{C}$ has duals, the horizontal functors have right adjoints by \cref{monfctdual,dualadjoint}; therefore by \cref{leftrightbc} one can equivalently check the right Beck-Chevalley condition, which reduces to the isomorphism
\[
f^{*}(A\tensor B^{d})\xto{\mu^{-1}} f^{*}A\tensor f^{*}B^{d}
\]
and we are done.
\endproof
\end{prop}

The dual notion and argument now gives

\begin{prop}[Right Projection Formula]\label{rproj} Let $f^{*}\colon \cat{C}\to\cat{D}$ be a monoidal functor, having a right adjoint $\ra{f}$. Then, there is a natural map 
\[
\rho\colon A\tensor\ra{f}B \to \ra{f}(f^{*}A\tensor B)
\]
If $\cat{C}$ has duals, then $\rho$ is an isomorphism. 
\end{prop}

One could  see that the projection morphisms satisfy several interesting properties such as compatibility with tensor product and composition, which can be understood as making the assignments $f\mapsto\lambda_{f}$, $f\mapsto\rho_{f}$ into ``monoidal functors''.   
However we will omit them for brevity and only consider the following, stating an invariance  of $\lambda$ and $\rho$  with respect to (op)lax squares, which will be needed later on. In particular, it implies stability under homotopy equivalences of functors, as one can easily check.

\begin{prop}\label{lambdamate}
Let
\[
\xymatrix{
\cat{A} \ar[d]_{h^{*}} \ar[r]^{g^{*}} & \cat{B} \ar[d]^{p^{*}}  \\
\cat{C} \ar[r]_{q^{*}}     & \cat{D} \ultwocell<\omit>{\phi}
}
\]
be a lax square of monoidal functors, commuting up to a monoidal transformation $\phi$.   
Assume  $\la{h}\lad h^{*}, \la{p}\lad p^{*}$.  
Then, the following diagram commutes
\[
\resizebox{\textwidth}{!}{
\xymatrix{
\la{p}(q^{*}A\tensor q^{*}h^{*}B) \ar[d]_{\la{p}(id\tensor\phi)} \ar[rr]^{\iso} &&   \la{p}q^{*}(A\tensor h^{*}B) \ar[rr]^{\lma{\phi}} && g^{*}\la{h}(A\tensor h^{*}B) \ar[d]^{g^{*}\lambda}  \\
\la{p}(q^{*}A\tensor p^{*}g^{*}B) \ar[r]_{\lambda} & \la{p}q^{*}A\tensor g^{*}B  \ar[rr]_{\lma{\phi}\tensor id} && g^{*}\la{h}A\tensor g^{*}B \ar[r]_{\iso} & g^{*}(\la{h}A\tensor B)
}
}
\]
\proof
Consider the diagrams
\[
\xymatrix@C+2pc{
\cat{A} \ar[d]_{h^{*}} \ar[r]^{-\tensor B} & \cat{A} \ar[d]_{h^{*}} \ar[r]^{g^{*}} & \cat{B} \ar[d]^{p^{*}} \\
\cat{C} \ar[d]_{id} \ar[r]_{-\tensor h^{*}B} & \cat{C} \ultwocell<\omit>{\mu} \ar[r]_{q^{*}} & \cat{D} \ultwocell<\omit>{\phi} \ar[d]^{id} \\
\cat{C} \ar[r]_{q^{*}} 
 \ar@{}[urr]^(.43){}="a"^(.57){}="b" \ar@{=>}_{\mu}  "a";"b"
& \cat{D} \ar[r]_{-\tensor q^{*}h^{*}B} 
&\cat{D} 
}
\]
and 
\[
\xymatrix@C+2pc{
\cat{A} \ar[r]^{-\tensor B} \ar[d]_{id} 
& \cat{A} \ar[r]^{g^{*}} &\cat{B} \ar[d]^{id} \\
\cat{A} \ar[d]_{h^{*}} \ar[r]^{g^{*}} 
 \ar@{}[urr]^(.43){}="a"^(.57){}="b" \ar@{=>}_{\mu}  "a";"b"
& \cat{B} \ar[d]_{p^{*}} \ar[r]^{-\tensor g^{*}B} & \cat{B} \ar[d]^{p^{*}} 
\\
\cat{C} \ar[d]_{id} \ar[r]_{q^{*}} & \cat{D} \ultwocell<\omit>{\phi} \ar[r]_{-\tensor p^{*}g^{*}B} & \cat{D}  \ultwocell<\omit>{\mu} \ar[d]^{id} \\
\cat{C} \ar[rr]_{q^{*}-\tensor q^{*}h^{*}B} 
 \ar@{}[urr]^(.43){}="a"^(.57){}="b" \ar@{=>}_{id\tensor\phi}  "a";"b"
& 
 & \cat{D} 
}
\]

whose left mates give respectively the upper and lower paths, and finally apply \cref{pasteleft}.
\endproof
\end{prop}

By reversing arrows and and using the pasting property for right mates we obtain

\begin{prop}\label{romate}
Let
\[
\xymatrix{
\cat{A} \ar[d]_{h^{*}} \ar[r]^{g^{*}} \drtwocell<\omit>{\pi} & \cat{B} \ar[d]^{p^{*}}  \\
\cat{C} \ar[r]_{q^{*}}     & \cat{D} 
}
\]
be an oplax square of monoidal functors, commuting up to a monoidal transformation $\pi$.   
Assume $h^{*}\lad\ra{h}, p^{*}\lad\ra{p}$.  
Then, the following diagram commutes
\[
\resizebox{\textwidth}{!}{
\xymatrix{
g^{*}(A\tensor\ra{h}B) \ar[r]^{\iso} \ar[d]_{g^{*}\rho} & g^{*}A\tensor g^{*}\ra{h}B \ar[rr]^{id\tensor\rma{\pi}} &&  g^{*}A\tensor \ra{p}q^{*}B \ar[r]^{\rho} & \ra{p}(p^{*}g^{*}A\tensor q^{*}B) \ar[d]^{\ra{p}(\pi\tensor id)} \\
g^{*}\ra{h}(h^{*}A\tensor B) \ar[rr]_{\rma{\pi}} && \ra{p}q^{*}(h^{*}A\tensor B) \ar[rr]_{\iso} && \ra{p}(q^{*}h^{*}A\tensor q^{*}B)
}
}
\]
\end{prop}

\subsection{Pre-Nakayama map} \label{prenakasubsec}

With the above results at hand, we can now introduce the Pre-Nakayama map $\ga{}$; once defined, we will study its behaviour with respect to homotopy commutative squares.

\begin{defin} \label{hgamma} 
Let $f^{*}\colon \cat{C}\to\cat{D}$ be a monoidal functor between monoidal categories with duals, having  both a right and a left adjoint. Define a natural map
\[
\hga{f}\colon\ra{f}(A\tensor B)\to \la{f} (A\tensor B)
\]
as 
\[
\resizebox{\textwidth}{!}{
\xymatrix{
\ra{f}(A\tensor B)\ar[rr]^-{\ra{f}(\etal{}\tensor id)} && \ra{f}(f^{*}\la{f}A\tensor B) \ar[r]^-{\lambda^{-1}\rho^{-1}} 
& \la{f}(A\tensor f^{*}\ra{f}B) \ar[rr]^-{\la{f} (id\tensor\epsr{})} &&\la{f}(A\tensor B)
}
}
\]
\end{defin}

\begin{defin}[Pre-Nakayama map] \label{ga}  
Let $f^{*}\colon \cat{C}\to\cat{D}$ be a monoidal functor between monoidal categories with duals, having  both a right and a left adjoint.  
Define a natural transformation, the \emph{Pre-Nakayama map},
\[
\ga{f}\colon\ra{f}\to \la{f}
\]
 as
\[
\ra{f}A\iso\ra{f}(A\tensor1)\xto{\hga{f}} \la{f}(A\tensor1)\iso\la{f}A
\]
\end{defin}

According to expectations, the map $\hga{}$, and consequently $\ga{}$, turns out to be well behaved w.r.t. homotopy commutative squares, by integrating \cref{lambdamate,romate} into 

\begin{lem}\label{hgamate}
Let
\[
\xymatrix@!0{
&  \cat{A} \ar[ld]_{g^{*}} \ar[rd]^{h^{*}} \\
\cat{B} \ar[rd]_{p^{*}}  &&   \cat{C}  \ar[ld]^{q^{*}}   \\
&  \cat{D}  \uutwocell<\omit>{\pi}
}
\]
be a diagram of monoidal functors between monoidal categories with duals, commuting up to a monoidal isomorphism $\pi$. Assume that
\begin{itemize}
\item the functors $h^{*},p^{*}$ have right and left adjoints
\item the square satisfies the left Beck-Chevalley condition relative to $\pi^{-1}$
\item the square satisfies the right Beck-Chevalley condition relative to $\pi$
\end{itemize}
Then, the following diagram commutes
\[
\xymatrix{
\ra{p}(q^{*}A\tensor q^{*}B)  \ar[r]^{\hga{p}}  &    \la{p}(q^{*}A\tensor q^{*}B) \ar[d]^{\iso} \\
\ra{p}q^{*}(A\tensor B) \ar[u]^{\iso}  &  \la{p}q^{*}(A \tensor B)  \ar[d]^{{\lma{(\pi^{-1})}}}  \\
g^{*}\ra{h}(A\tensor B) \ar[u]^{{\rma{\pi}}} \ar[r]_{g^{*}\hga{h}} & g^{*}\la{h}(A\tensor B) 
}
\]
\proof
The square arises by  gluing
\[
\xymatrix{
& \ra{p}(p^{*}\la{p}q^{*}A\tensor q^{*}B) \ar[r]^{\rho^{-1}} 
 \ar[d]^{\ra{p}(p^{*}{\lma{(\pi^{-1})}}\tensor id)}
 &  \la{p}q^{*}A\tensor\ra{p}q^{*}B  \\
\ra{p}(q^{*}A\tensor q^{*}B)  \ar[rd]_{\ra{p}(q^{*}\etal{h}\tensor id)}  \ar[ru]^{\ra{p}(\etal{p}\tensor id)}  & \ra{p}(p^{*}g^{*}\la{h}A\tensor q^{*}B)  
\ar[r]_{\rho^{-1}} & g^{*}\la{h}A\tensor \ra{p}q^{*}B \ar[u]_{{\lma{(\pi^{-1})}}^{-1}\tensor id}    \\
&  \ra{p}(q^{*}h^{*}\la{h}A\tensor q^{*}B)   \ar[u]_{\ra{p}(\pi^{-1}\tensor id)}  & \\
\ra{p}q^{*}(A\tensor B) \ar[uu]^{\iso} 
&    \ra{p}q^{*}(h^{*}\la{h}A\tensor B)  \ar[u]^{\iso} 
\ar@{}[ru]|{\cref{romate}}
& 
g^{*}\la{h}A\tensor g^{*}\ra{h}B    \ar[uu]_{id\tensor {\rma{\pi}}}
\\
g^{*}\ra{h}(A\tensor B) \ar[u]^{{\rma{\pi}}} \ar[r]_{g^{*}\ra{h}(\etal{h}\tensor id)}  & g^{*}\ra{h}(h^{*}\la{h}A\tensor B) \ar[u]^{{\rma{\pi}}} \ar[r]_{g^{*}\rho^{-1}} & g^{*}(\la{h}A\tensor \ra{h}B)    \ar[u]_{\iso}
}
\]

with its specular diagram obtained via $\lambda^{-1}$. 
\cref{lambdamate,romate} and the Beck-Chevalley condition hypothesis on $\pi$ and $\pi^{-1}$, together with a short computation involving the definition of mates and the axioms for monoidal structures, conclude the proof.
%
%
%
%
%
\endproof
\end{lem}

By applying the structure isomorphisms for monoidal categories and functors, one easily extends the above lemma to $\ga{}$

\begin{prop} \label{gamate}
Let
\[
\xymatrix@!0{
&  \cat{A} \ar[ld]_{g^{*}} \ar[rd]^{h^{*}} \\
\cat{B} \ar[rd]_{p^{*}}  &&   \cat{C}  \ar[ld]^{q^{*}}   \\
&  \cat{D}  \uutwocell<\omit>{\pi}
}
\]
be a diagram of monoidal functors between monoidal categories with duals, commuting up to a monoidal isomorphism $\pi$. Assume that
\begin{itemize}
\item the functors $h^{*},p^{*}$ have right and left adjoints
\item the square satisfies the left Beck-Chevalley condition relative to $\pi^{-1}$
\item the square satisfies the right Beck-Chevalley condition relative to $\pi$
\end{itemize}
Then, the following diagram commutes
\[
\xymatrix{
\ra{p}q^{*}A  \ar[r]^{\ga{p}}  &    \la{p}q^{*}A \ar[d]^{{\lma{(\pi^{-1})}}} \\
g^{*}\ra{h}A \ar[u]^{\rma{\pi}} \ar[r]_{g^{*}\ga{h}} & g^{*}\la{h}A
}
\]
%
%
%
\end{prop}

In particular for $g^{*},q^{*}$ identities, it is easy to see that the above exhibit $\hga{}$ and $\ga{}$ as stable, in a suitable sense, under equivalences of functors. 
Notice also that, due to \cref{montraiso}, one doesn't need to require $\pi$ to be an isomorphism; we preferred however to emphasize its invertibility, which would anyway be needed in a non-dualizable setting.\\

We shall now prove a technical result, which will be the keystone for the functoriality of the quantization functor $\sam$, and states the equivalence of the two possible ``paths'' one can walk along a homotopy commutative square.

\begin{thm}  \label{pullpushga}
Let
\[
\xymatrix@!0{
&  \cat{A} \ar[ld]_{g^{*}} \ar[rd]^{h^{*}} \\
\cat{B} \ar[rd]_{p^{*}}  &&   \cat{C}  \ar[ld]^{q^{*}}   \\
&  \cat{D}  \uutwocell<\omit>{\pi}
}
\]
be a diagram of monoidal functors between monoidal categories with duals, commuting up to a monoidal isomorphism $\pi$. Assume that
\begin{itemize}
\item the functors $h^{*},p^{*}$ have right and left adjoints
\item the functors $g^{*},q^{*}$ have left adjoints
\item the square satisfies the left Beck-Chevalley condition relative to $\pi^{-1}$
\item the square satisfies the right Beck-Chevalley condition relative to $\pi$
\end{itemize}
Let $\bar{\pi}\colon \la{g}\la{p}\to\la{h}\la{q}$ be the natural isomorphism induced by $\pi$.  
Then, the following diagram commutes
\[
\resizebox{\textwidth}{!}{
\xymatrix@C+1pc{
\la{g}g^{*} \ar[rrd]_{\etar{h}} \ar[r]^{\la{g}\etar{p}} &  \la{g}\ra{p}p^{*}g^{*} \ar[r]^{\la{g}\ra{p}\pi}  & \la{g}\ra{p}q^{*}h^{*}  \ar[r]^{\la{g}\ga{p}} &  \la{g}\la{p}q^{*}h^{*} \ar[r]^{\bar{\pi}} &  \la{h}\la{q}q^{*}h^{*} \ar[r]^{\la{h}\epsl{q}} &  \la{h}h^{*} \\
&& \ra{h}h^{*}\la{g}g^{*} \ar[r]_{\ga{h}} &  \la{h}h^{*}\la{g}g^{*} \ar[rru]_{\la{h}h^{*}\epsl{g}} &&
}
}
\]
\proof
Noting that $\bar{\pi}$ is the composition
\[
\la{g}\la{p}\xto{\la{g}\la{p}\etal{q}}\la{g}\la{p}q^{*}\la{q}
\xto{\la{g}\lma{(\pi^{-1})}}\la{g}g^{*}\la{h}\la{q}\xto{\epsl{g}}\la{h}\la{q}
\]
the proof follows from \cref{gamate}, naturality of maps and a short computation involving the definition of mates.
\endproof
\end{thm}

As an aside observation, we compare a map given in \cite{fhm} with the above $\ga{}$. Again, we ask for a category with duals, even if the result holds (as for the previous ones) for functors satisfying the projection formulas.
\begin{defin} \label{mayga} 
Let $f^{*}\colon \cat{C}\to\cat{D}$ be a monoidal functor between categories with duals, having  both a right and a left adjoint, and assume we have a map $\chi\colon\ra{f}1\to\la{f}1$. Define a natural transformation
\[
\mu_{\chi} \colon \ra{f}\to\la{f}
\]
as
\begin{align*}
\ra{f}X\iso\ra{f}X\tensor1 & \xto{id\tensor\etar{}}\ra{f}X\tensor\ra{f}f^{*}1 \xto{id\tensor\chi} \ra{f}X\tensor\la{f}f^{*}1\iso \la{f}f^{*}1\tensor\ra{f}X \\
& \xto{\lambda^{-1}}\la{f}(f^{*}1\tensor f^{*}\ra{f}X)\iso\la{f}f^{*}\ra{f}X\xto{\la{f}\epsr{}}\la{f}X
\end{align*}
\end{defin}

\begin{prop} \label{gamay}
Let $f^{*}\colon \cat{C}\to\cat{D}$ be a monoidal functor between categories with duals, having  both a right and a left adjoint. Then
\[
\mu_{\ga{f}(1)}=\ga{f}
\]
\proof
The proposition follows from the properties of $\ga{}$ and of the projection maps.
\endproof
\end{prop}

\subsection{Tensoring Kan extensions and pre-Nakayama maps} \label{prenakakantensorsubsec}\label{nakakantensor} 

In the following we will focus on external tensor product of functors and how this is respected both by left and right Kan extensions, in the sense that we have isomorphisms
\begin{align*}
\la{\mu} &\colon \la{(f\times g)}(A\etensor B) \to \la{f}A\etensor \la{g}B \\
\ra{\mu} &\colon \ra{(f\times g)}(A\etensor B) \to \ra{f}A\etensor \ra{g}B 
\end{align*}
relating the Kan extensions along a cartesian product of maps and the tensor product of the single Kan extensions. We will then see in \cref{multiplyga} that the pre-Nakayama map agrees with such structure.\\

Let us first recall the external tensor product of functors, which might be understood as the categorical interpretation of external tensor product of group representations, and can be defined on presheaves valued in a monoidal category as follows

\begin{defin} \label{exttensor}
Let $\cat{C}$ be a monoidal category, and $A\in\cat{C}^{X},B\in\cat{C}^{Y}$ functors. The \emph{external tensor product} of $A$ and $B$ is the functor $A\etensor B\in\cat{C}^{X\times Y}$ 
\[
A\etensor B \defeq 
 p_{X}^{*}A\tensor p_{Y}^{*}B
\]
where $p_{X},p_{Y}$ are the obvious projections.  

One has natural ``associativity'' and ``twist'' isomorphisms
\begin{align*}
\bar{\alpha}\colon & (A\etensor B)\etensor C \to A\etensor (B\etensor C) \\
\bar{\tau}\colon & A\etensor B\to B\etensor A
\end{align*}
given respectively by 

\[
\resizebox{\textwidth}{!}{
\xymatrix@C+2pc@R+1pc{
\cat{C}^{X}\times\cat{C}^{Y}\times\cat{C}^{Z} \ar[r]^-{id \times p_{Y}^{*} \times p_{Z}^{*}} \ar[d]_{p_{X}^{*} \times p_{Y}^{*} \times id}  &   \cat{C}^{X}\times\cat{C}^{Y\times Z}\times \cat{C}^{Y\times Z} \ar[r]^-{id \times \tensor}  \ar[d]|{p_{X}^{*} \times p_{Y\times Z}^{*} \times p_{Y\times Z}^{*}}    &     \cat{C}^{X}\times\cat{C}^{Y \times Z} \ar[d]^{p_{X}^{*} \times p_{Y\times Z}^{*}} \\
\cat{C}^{X\times Y} \times \cat{C}^{X\times Y} \times \cat{C}^{Z} \ar[r]_-{p_{X\times Y}^{*} \times p_{X\times Y}^{*} \times p_{Z}^{*}} \ar[d]_{\tensor \times id}  & \cat{C}^{X\times Y\times Z} \times \cat{C}^{X\times Y\times Z} \times \cat{C}^{X\times Y\times Z}
 \ar@{}[ur]^(.43){}="a"^(.57){}="b" \ar@{=>}_{\mu^{-1}}  "a";"b"
 \ar[d]|{\tensor \times id} \ar[r]^-{id \times \tensor}   & \cat{C}^{X\times Y \times Z}\times \cat{C}^{X\times Y \times Z}  \ar[d]^{\tensor}  
\\
\cat{C}^{X\times Y} \times \cat{C}^{Z} \ar[r]_-{p_{X\times Y}^{*} \times p_{Z}^{*}} 
 \ar@{}[ur]^(.43){}="a"^(.57){}="b" \ar@{=>}_{\mu}  "a";"b"
 &  \cat{C}^{X\times Y\times Z} \times \cat{C}^{X\times Y\times Z} \ar[r]_-{\tensor}
 \ar@{}[ur]^(.43){}="a"^(.57){}="b" \ar@{=>}_{\alpha}  "a";"b"
& \cat{C}^{X\times Y\times Z}
}
}
\]
and
\[
\xymatrix{
\cat{C}^{X}\times\cat{C}^{Y} \ar[r]^-{p_{X}^{*}\times p_{Y}^{*}} \ar[d]_{id} & \cat{C}^{X\times Y} \times \cat{C}^{X\times Y} \ar[r]^-{\tensor}\ar[d]_{\iso} & \cat{C}^{X\times Y}\ar[d] \\
\cat{C}^{X}\times\cat{C}^{Y} \ar[r]_-{p_{X}^{*}\times p_{Y}^{*}} \ar[d]_{\iso} & \cat{C}^{Y\times X} \times \cat{C}^{Y\times X} \ar[r]^-{\tensor}\ar[d]_{id} & \cat{C}^{Y\times X} \ultwocell<\omit>{\mu} \ar[d]^{id} \\
\cat{C}^{Y}\times\cat{C}^{X} \ar[r]_-{p_{Y}^{*}\times p_{X}^{*}}  & \cat{C}^{Y\times X} \times \cat{C}^{Y\times X} \ar[r]_-{\tensor} & \cat{C}^{Y\times X}  \ultwocell<\omit>{\tau}
}
\] 
where we considered the obvious twisting isomorphisms for the cartesian product.
\end{defin}

From now $\cat{C}$ will be a fixed monoidal category with duals, closed under colimits and limits indexed by groupoids. This will ensure existence of both left and right Kan extensions along any map $f$ between groupoids, which we denote by $\la{f}$ and $\ra{f}$; groupoids will be denoted by capital letters $M,N,\dots,X,Y,\dots$ 
Also, recall from \cref{fundualcat} that the functor categories $\cat{C}^{X}$ are symmetric monoidal with duals. 

The primary goal is to define isomorphisms $\la{\mu},\ra{\mu}$. One should think of the following as making the Kan extension constructions into symmetric monoidal functors. 

\begin{prop} \label{multiply}
Given maps $f\colon M\to X$ and $g\colon N\to Y$, for any  $A\in \cat{C}^{M},B\in\cat{C}^{N}$ we have natural isomorphisms
\begin{align*}
\la{\mu} &\colon \la{(f\times g)}(A\etensor B) \to \la{f}A\etensor \la{g}B \\
\ra{\mu} &\colon \ra{(f\times g)}(A\etensor B) \to \ra{f}A\etensor \ra{g}B 
\end{align*}
\proof
We prove the existence of the isomorphism $\la{\mu}$, being the other dual. Consider the left mate resulting from the following diagram
\[
\xymatrix@C+2pc @R+1pc{
\cat{C}^{X}\times \cat{C}^{Y} \ar[d]_{f^{*}\times id} \ar[r]^{p_{X}^{*}\times p_{Y}^{*}}  \ar@/^2pc/[rr]^{\boxtimes}  & \cat{C}^{X\times Y} \times \cat{C}^{X\times Y} \ar[r]^{\tensor} \ar[d]|{(f\times id)^{*}\times id} & \cat{C}^{X\times Y} \ar[d]^{(f\times id)^{*}} \\
\cat{C}^{M} \times \cat{C}^{Y} \ar[d]_{id} \ar[r]^{p_{M}^{*}\times p_{Y}^{*}}   &  \cat{C}^{M\times Y} \times  \cat{C}^{X\times Y} \ar[r]^-{\tensor (f\times id)^{*}}  &   \cat{C}^{M\times Y} \ar[d]^{id} \ultwocell<\omit>{\mu} \\
\cat{C}^{M} \times \cat{C}^{Y} \ar[d]_{id\times g^{*}} \ar[r]^{p_{M}^{*}\times p_{Y}^{*}}   &  \cat{C}^{M\times Y} \times  \cat{C}^{M\times Y} \ar[r]^-{\tensor} \ar[d]|{id\times (id\times g)^{*}}   &   \cat{C}^{M\times Y} \ar[d]^{(id\times g)^{*}}\\
 \cat{C}^{M}\times  \cat{C}^{N} \ar[r]_{p_{M}^{*}\times p_{N}^{*}} \ar@/_2pc/[rr]_{\boxtimes}  &  \cat{C}^{M\times Y} \times \cat{C}^{M\times N} \ar[r]_-{(id\times g)^{*}\tensor}  \ultwocell<\omit>{id} & \cat{C}^{M\times N} \ultwocell<\omit>{\mu}\\
 &&
} 
\]
The left sided squares give rise to isos on each factor thanks to \cref{grpBC}. The right sided squares reduce, pointwise, to the maps $\lambda$ of \cref{lproj}.
\endproof
\end{prop}

A simple argument then proves associativity of $\la{\mu}$ and $\ra{\mu}$

\begin{prop} \label{multiplyassociativity}
Given maps $M\xto{f}X$, $N\xto{g}Y$, $O\xto{h}Z$, for any ${A\in \cat{C}^{M}},{B\in\cat{C}^{N}},{C\in\cat{C}^{O}}$ the following diagrams commute

\[
\resizebox{\textwidth}{!}{
\xymatrix@C+1pc{
\la{(f\times g\times h)}((A\etensor B)\etensor C) \ar[d]_{\iso} \ar[r]^{\la{\mu}} &
\la{(f\times g)}(A\etensor B) \etensor \la{h} C \ar[r]^{\la{\mu} \etensor id} & (\la{f}A\etensor\la{g}B) \etensor \la{h}C
  \ar[d]^{\iso} \\
\la{(f\times g\times h)}( A\etensor (B \etensor C)) \ar[r]_{\la{\mu}} & 
\la{f}A\etensor\la{(g\times h)}(B\etensor C) \ar[r]_{id\etensor \la{\mu}} & \la{f}A\etensor (\la{g}B\etensor \la{h}C) 
}
}
\]

\[
\resizebox{\textwidth}{!}{
\xymatrix@C+1pc{
\ra{(f\times g\times h)}(A\etensor(B\etensor C)) \ar[d]_{\iso} \ar[r]^{\ra{\mu}} & \ra{f}A\etensor\ra{(g\times h)}(B\etensor C) \ar[r]^{id\etensor \ra{\mu}} & \ra{f}A\etensor (\ra{g}B\etensor \ra{h}C) \ar[d]^{\iso} \\
\ra{(f\times g\times h)}((A\etensor B)\etensor C) \ar[r]_{\ra{\mu}} & \ra{(f\times g)}(A\etensor B) \etensor \ra{h} C \ar[r]_{\ra{\mu} \etensor id} & (\ra{f}A\etensor\ra{g}B) \etensor \ra{h}C 
}
}
\]

\proof
Being the two diagrams dual to each other, we give just the proof for the first. 
The upper and lower paths are, respectively, the left mates of  
\[
\xymatrix@C+1.5pc{
\cat{C}^{X}\times\cat{C}^{Y}\times\cat{C}^{Z} \ar[d]|{id} \ar[r]^{id \times \etensor} & \cat{C}^{X}\times \cat{C}^{Y\times Z} \ar[r]^{\etensor} & \cat{C}^{X\times Y \times Z} \ar[d]|{id^{*}} \\
\cat{C}^{X}\times\cat{C}^{Y}\times\cat{C}^{Z} \ar[d]|{f^{*}\times id \times id} \ar[r]^{\etensor\times id} & \cat{C}^{X\times Y}\times \cat{C}^{Z} \ar[d]|{(f\times id)^{*}\times id} \ar[r]^{\etensor} & \cat{C}^{X\times Y \times Z} \ar[d]|{(f\times id\times id)^{*}} \ulltwocell<\omit>{\bar{\alpha}} \\
\cat{C}^{M}\times\cat{C}^{Y}\times\cat{C}^{Z} \ar[d]|{id\times g^{*}\times id} \ar[r]^{\etensor\times id} & \cat{C}^{M\times Y}\times \cat{C}^{Z} \ar[r]^{\etensor} \ar[d]|{(id\times g)^{*}\times id} \ultwocell<\omit>{} & \cat{C}^{M\times Y \times Z} \ar[d]|{(id\times g \times id)^{*}} \ultwocell<\omit>{} \\
\cat{C}^{M}\times\cat{C}^{N}\times\cat{C}^{Z} \ar[d]|{id\times id\times h^{*}} \ar[r]^{\etensor\times id} & \cat{C}^{M\times N}\times \cat{C}^{Z} \ar[d]|{(id\times id) \times h^{*}} \ar[r]^{\etensor} \ultwocell<\omit>{} & \cat{C}^{M\times N \times Z} \ar[d]|{(id\times id\times h)^{*}} \ultwocell<\omit>{} \\
\cat{C}^{M}\times\cat{C}^{N}\times\cat{C}^{O} \ar[r]^{\etensor\times id} & \cat{C}^{M\times N}\times \cat{C}^{O} \ar[r]^{\etensor} \ultwocell<\omit>{} & \cat{C}^{M\times N \times O} \ultwocell<\omit>{}\\
}
\]

and

\[
\xymatrix@C+1.5pc{
\cat{C}^{X}\times\cat{C}^{Y}\times\cat{C}^{Z} \ar[d]|{f^{*}\times id \times id} \ar[r]^{id \times \etensor} & \cat{C}^{X}\times \cat{C}^{Y\times Z} \ar[d]|{f^{*}\times(id\times id)} \ar[r]^{\etensor} & \cat{C}^{X\times Y \times Z} \ar[d]|{(f\times id\times id)^{*}} \\ 
\cat{C}^{M}\times\cat{C}^{Y}\times\cat{C}^{Z} \ar[d]|{id\times g^{*}\times id} \ar[r]^{id \times \etensor} & \cat{C}^{M}\times \cat{C}^{Y\times Z} \ar[r]^{\etensor} \ar[d]|{id\times (g\times id)^{*}} \ultwocell<\omit>{} & \cat{C}^{M\times Y \times Z} \ar[d]|{(id\times g \times id)^{*}} \ultwocell<\omit>{} \\
\cat{C}^{M}\times\cat{C}^{N}\times\cat{C}^{Z} \ar[d]|{id\times id\times h^{*}} \ar[r]^{id \times \etensor} & \cat{C}^{M}\times \cat{C}^{N\times Z} \ar[d]|{id\times (id \times h)^{*}} \ar[r]^{\etensor} \ultwocell<\omit>{} & \cat{C}^{M\times N \times Z} \ar[d]|{(id\times id\times h)^{*}} \ultwocell<\omit>{} \\
\cat{C}^{M}\times\cat{C}^{N}\times\cat{C}^{O} \ar[d]|{id} \ar[r]^{id \times \etensor} & \cat{C}^{M}\times \cat{C}^{N\times O} \ar[r]^{\etensor} \ultwocell<\omit>{} & \cat{C}^{M\times N \times O} \ultwocell<\omit>{} \ar[d]|{id}\\
\cat{C}^{M}\times\cat{C}^{N}\times\cat{C}^{O}  \ar[r]^{\etensor\times id} & \cat{C}^{M\times N}\times \cat{C}^{O} \ar[r]^{\etensor} & \cat{C}^{M\times N \times O}  \ulltwocell<\omit>{\bar{\alpha}}
}
\]
where we pasted the diagrams appearing in \cref{exttensor,multiply}. They are easily seen to be equivalent by definition of $\etensor$ and the transformations involved, and monoidality of the restriction functors. \cref{pasteleft} then concludes the proof.
\endproof
\end{prop}

Again, via the pasting lemmas, it is not difficult to check the symmetry axiom.

\begin{prop} \label{multiplycommute}
Given maps $M\xto{f}X$, $N\xto{g}Y$, for any $A\in \cat{C}^{M},B\in\cat{C}^{N}$ the following diagrams commute

\[
\xymatrix{
\la{(f\times g)}(A\etensor B) \ar[d]_{\iso}  \ar[r]^-{\la{\mu}} & \la{f}A\etensor\la{g}B \ar[d]^{\iso} \\
\la{(g\times f)}(B\etensor A) \ar[r]_-{\la{\mu}} & \la{g}B\etensor\la{f}A
}
\]

\[
\xymatrix{
\ra{(f\times g)}(A\etensor B) \ar[d]_{\iso}  \ar[r]^-{\ra{\mu}} & \ra{f}A\etensor\ra{g}B \ar[d]^{\iso} \\
\ra{(g\times f)}(B\etensor A) \ar[r]_-{\ra{\mu}} & \ra{g}B\etensor\ra{f}A
}
\]
\end{prop}

The same arguments  also show that $\la{\mu},\ra{\mu}$ behave functorially

\begin{prop}\label{multiplycomp}
Let $P\xto{p}M\xto{f}X$ and $Q\xto{q}N\xto{g}Y$ be pairs of composable maps. Then for any ${A\in \cat{C}^{P}},{B\in\cat{C}^{Q}}$, the diagrams

\[
\resizebox{\textwidth}{!}{
\xymatrix@C+2.5pc @R+1pc{
\la{(f\times g)}\la{(p\times q)}(A\etensor B) \ar[r]^{\la{(f\times g)}\la{\mu}}   \ar@/_2pc/[rr]_{\la{\mu}} 
& \la{(f\times g)}(\la{p}A\etensor \la{q}B) \ar[r]^{\la{\mu}} &    \la{f}\la{p}A\etensor \la{g}\la{q}B
}
}
\]

\[
\resizebox{\textwidth}{!}{
\xymatrix@C+2.5pc @R+1pc{
\ra{(f\times g)}\ra{(p\times q)}(A\etensor B) \ar[r]^{\ra{(f\times g)}\ra{\mu}}   \ar@/_2pc/[rr]_{\ra{\mu}} 
& \ra{(f\times g)}(\ra{p}A\etensor \ra{q}B) \ar[r]^{\ra{\mu}} &
  \ra{f}\ra{p}A\etensor \ra{g}\ra{q}B
}
}
\]
commute.
\proof
The first composition results by pasting the left mates of the inner squares in the following diagram, while the lower map is the left mate of the outer diagram:
\[
\xymatrix@R+1.5pc{
\cat{C}^{X}\times\cat{C}^{Y} \ar[d]_{f^{*}\times g^{*}} \ar[r]^{\etensor} & \cat{C}^{X\times Y} \ar[d]^{(f\times g)^{*}}   \\
\cat{C}^{M}\times\cat{C}^{N} \ar[d]_{p^{*}\times q^{*}} \ar[r]^{\etensor} & \cat{C}^{M\times N} \ar[d]^{(p\times q)^{*}} \ultwocell<\omit>{} \\
\cat{C}^{P}\times\cat{C}^{Q} \ar[r]_{\etensor} & \cat{C}^{P\times Q} \ultwocell<\omit>{}
}
\]
\endproof
\end{prop}

Finally, as one would expect, stability under homotopy equivalences is an immediate consequence of the Beck-Chevalley formalism.

\begin{prop}\label{multiplyhomo}
Let
\[
\xymatrix{
M\rtwocell^{f}_{f'}{\phi} & X
}
\qquad
\xymatrix{
N\rtwocell^{g}_{g'}{\psi} & X
}
\]
be natural transformations of maps between groupoids. Then for any ${A\in \cat{C}^{M}},{B\in\cat{C}^{N}}$, the following diagrams commute

\[
\xymatrix{
\la{(f'\times g')}(A\etensor B) \ar[d]_{\lma{(\phi\times\psi)}} \ar[r]^-{\la{\mu}}  &   \la{f'}A\etensor\la{g'}B  \ar[d]^{\lma{\phi}\etensor\lma{\psi}}\\
\la{(f\times g)}(A\etensor B) \ar[r]^-{\la{\mu}}  &   \la{f}A\etensor\la{g}B 
}
\]

\[
\xymatrix{
\ra{(f'\times g')}(A\etensor B)\ar[d]_{\rma{(\phi\times\psi)}} \ar[r]^-{\ra{\mu}}  &   \ra{f'}A\etensor\ra{g'}B  \ar[d]^{\rma{\phi}\etensor\rma{\psi}} \\
\ra{(f\times g)}(A\etensor B) \ar[r]^-{\ra{\mu}}  &   \ra{f}A\etensor\ra{g}B 
}
\]
\proof
The upper and lower paths in the first square are given by the left mates of the following diagrams
\[
\xymatrix{
\cat{C}^{X}\times\cat{C}^{Y} \ar[r]^{id}\ar[d]_{f^{*}\times g^{*}} & \cat{C}^{X}\times\cat{C}^{Y} \ar[d]|{f'^{*}\times g'^{*}} \ar[r]^{\etensor}  & \cat{C}^{X\times Y} \ar[d]^{(f'\times g')^{*}} \\
\cat{C}^{M}\times\cat{C}^{N} \ar[r]_{id} 
\ar@{}[ur]^(.38){}="a"^(.62){}="b" \ar@{=>}^{\phi^{*}\times \psi^{*}}  "a";"b"
& \cat{C}^{M}\times\cat{C}^{N}  \ar[r]_{\etensor} 
\ar@{}[ur]^(.38){}="a"^(.62){}="b" \ar@{=>} "a";"b"
& \cat{C}^{M\times N} 
}
\]
\[
\xymatrix{
\cat{C}^{X}\times\cat{C}^{Y} \ar[d]_{f^{*}\times g^{*}} \ar[r]^{\etensor}  & \cat{C}^{X\times Y} \ar[d]|{(f\times g)^{*}} \ar[r]^{id} &  \cat{C}^{X\times Y} \ar[d]^{(f'\times g')^{*}}\\
\cat{C}^{M}\times\cat{C}^{N}  
\ar@{}[ur]^(.38){}="a"^(.62){}="b" \ar@{=>}  "a";"b"
\ar[r]_{\etensor}   & \cat{C}^{M\times N}  
 \ar@{}[ur]^(.38){}="a"^(.62){}="b" \ar@{=>}_{(\phi\times \psi)^{*}}  "a";"b"
 \ar[r]_{id} & \cat{C}^{M\times N}  
}
\]
which are easily seen to be equivalent, for $\phi$ and $\psi$ induce monoidal transformations. 

Commutativity of the second square is proved similarly.
\endproof
\end{prop}

In view of \cref{thmsum,thmprod} we also provide the following

\begin{prop}\label{multiplyetaeps}
Let $f\colon M\to X$ and $g\colon N\to Y$ be maps. Then for any ${A\in \cat{C}^{X}},{B\in\cat{C}^{Y}}$, the diagrams

\[
\resizebox{\textwidth}{!}{
\xymatrix@C+1.5pc{
\la{(f\times g)}(f\times g)^{*}(A\etensor B)\ar[r]^{\iso} \ar@/_2pc/[rrr]_{\eps}  & \la{(f\times g)}(f^{*}A\etensor g^{*}B) \ar[r]^{\la{\mu}}   
& \la{f}f^{*}A\etensor \la{g}g^{*}B \ar[r]^-{\eps\etensor\eps}   &  A\etensor B
}
}
\]

\[
\resizebox{\textwidth}{!}{
\xymatrix@C+1.5pc{
A\etensor B \ar[r]^-{\eta}  \ar@/_2pc/[rrr]_{\eta\etensor\eta}  
& \ra{(f\times g)}(f\times g)^{*}(A\etensor B)\ar[r]^{\iso} & \ra{(f\times g)}(f^{*}A\etensor g^{*}B) \ar[r]^{\ra{\mu}} & \ra{f}f^{*}A\etensor \ra{g}g^{*}B   
%
}
}
\]
commute.
\proof
A simple computation shows that the diagram
\[
\resizebox{\textwidth}{!}{
\xymatrix@C+2.5pc @R+1pc{
\cat{C}^{X}\times \cat{C}^{Y} \ar[d]_{id} \ar[r]^{id}    & \cat{C}^{X}\times \cat{C}^{Y} \ar[d]_{f^{*}\times g^{*}} \ar[r]^{p_{X}^{*}\times p_{Y}^{*}}  & \cat{C}^{X\times Y} \times \cat{C}^{X\times Y} \ar[r]^{\tensor} 
& \cat{C}^{X\times Y} \ar[d]^{(f\times g)^{*}} \\
\cat{C}^{X}\times \cat{C}^{Y} \ar[d]_{id} \ar[r]^{f^{*}\times g^{*}} & \cat{C}^{M}\times  \cat{C}^{N}
\ar@{}[urr]^(.46){}="a"^(.54){}="b" \ar@{=>}^{\mu}  "a";"b"
 \ar[r]^{p_{M}^{*}\times p_{N}^{*}}   &  \cat{C}^{M\times N} \times \cat{C}^{M\times N} \ar[d]^{id}  \ar[r]_{\tensor}   & \cat{C}^{M\times N} \ar[d]^{id} 
 \\
\cat{C}^{X}\times \cat{C}^{Y} \ar[d]_{id}  \ar[r]^{p_{X}^{*}\times p_{Y}^{*}}  &  \cat{C}^{X\times Y} \times \cat{C}^{X\times Y} \ar[d]^{id}  \ar[r]^{(f\times g)^{*}\times(f\times g)^{*}}   &   \cat{C}^{M\times N} \times \cat{C}^{M\times N} \ar[r]^{\tensor} &  \cat{C}^{M\times N} \ar[d]^{id} \\
\cat{C}^{X}\times \cat{C}^{Y}  \ar[r]_{p_{X}^{*}\times p_{Y}^{*}}  & \cat{C}^{X\times Y} \times \cat{C}^{X\times Y} 
\ar@{}[urr]^(.46){}="a"^(.54){}="b" \ar@{=>}^{\mu^{-1}}  "a";"b"
\ar[r]_{\tensor}  & \cat{C}^{X\times Y} \ar[r]_{(f\times g)^{*}} & \cat{C}^{M\times N} 
} 
}
\]
whose left mate gives the first composition, 
is equivalent to the diagram
\[
\xymatrix{
\cat{C}^{X\times Y} \ar[r]^{id} \ar[d]_{id} & \cat{C}^{X\times Y} \ar[d]^{(f\times g)^{*}} \\
\cat{C}^{X\times Y} \ar[r]_{(f\times g)^{*}}  & \cat{C}^{M\times N}
}
\]
giving the counit $\eps$.

The same argument proves the other statement.
\endproof
\end{prop}

We can now get back to our pre-Nakayama map, and conclude the section with a compatibility theorem

\begin{thm}\label{multiplyga}  
Let $f\colon M\to X$ and $g\colon N\to Y$ be maps. Then for any ${A\in \cat{C}^{M}},{B\in\cat{C}^{N}}$ 
\[
\xymatrix{
\ra{(f\times g)}(A\etensor B) \ar[d]_{\ga{f\times g}}  \ar[r]^{\ra{\mu}} &  \ra{f}A\etensor \ra{g}B \ar[d]^{\ga{f}\etensor\ga{g}} \\
\la{(f\times g)}(A\etensor B) \ar[r]_{\la{\mu}} & \la{f}A\etensor \la{g}B
}
\]
commutes.
\proof
By definition, $\la{\mu}$ and $\ra{\mu}$ are the left and right mates for the diagram 
\[
\xymatrix{
\cat{C}^{X}\times\cat{C}^{Y} \ar[d]_{f^{*}\times g^{*}} \ar[r]^{\etensor} & \cat{C}^{X\times Y} \ar[d]^{(f\times g)^{*}}   \\
\cat{C}^{M}\times\cat{C}^{N}  \ar[r]_{\etensor} & \cat{C}^{M\times N}  \ultwocell<\omit>{} \\
}
\]
Since all functors and transformations involved are monoidal, and since the pre-Nakayama map relative to $f^{*}\times g^{*}$ is $\ga{f}\times\ga{g}$, we can apply \cref{gamate}.
\endproof
\end{thm}

\section{Nakayama Categories} \label{nakacontextsec}\label{nakasubsec} \label{secnaka}

As already mentioned in the introduction and shown below, the assignment $f\mapsto\ga{f}$ fails to be functorial in general. The latter property will be essential to the construction of the quantization functors, in order to basically reduce our problem to \cref{pullpushga}. 

The next pages are therefore devoted to correct this issue, by building a new map $\na{}$, and to the definition of what we call a ``Nakayama category''; this can be understood as an abstract version of finite vector spaces over a zero-characteristic field, and will encode the hypotheses needed to build the functor $\sam$.\\

We give first the following, which allows us to reduce the assumptions on our category $\cat{C}$

\begin{prop}\label{rightiffleft}
Let $\cat{C}$ and $\cat{D}$ be monoidal categories with duals, and ${f^{*}\colon \cat{C}\to\cat{D}}$ a monoidal functor. Then, $f^{*}$ has a left adjoint if and only if it has a right adjoint.
\proof
Assume that $f^{*}$ has a left adjoint, then we have a chain of isomorphisms
\[
\Hom{\cat{D}}{f^{*}B}{A}\iso \Hom{\cat{D}}{A^{d}}{(f^{*}B)^{d}}\iso  \Hom{\cat{D}}{A^{d}}{f^{*}B^{d}} \iso   \Hom{\cat{C}}{\la{f}A^{d}}{B^{d}} \iso \Hom{\cat{C}}{B}{(\la{f}A^{d})^{d}}
\]
expressing $f^{*}$ as a left adjoint to the functor
\[
A\mapsto (\la{f}A^{d})^{d}
\]
The other direction of the equivalence is proved similarly.
\endproof
\end{prop}

From now on we will be interested in categories $\cat{C}^{X}$ of functors from essentially finite groupoids, by which we mean

\begin{defin}  
A groupoid $X$ is said \emph{essentially finite} if the set $\pi_{0}X$ of isomorphism classes of objects of $X$ is finite, and for each pair of objects $x,y\in X$ the set $\Hom{X}{x}{y}$ is finite.
\end{defin}

It can be  seen that the  Quillen equivalence between homotopy $1$-types and groupoids (see for example \cite[Section 2.2]{jt})  restricts to a Quillen equivalence between the categories of finite homotopy $1$-types and of essentially finite groupoids. We will denote both categories by $\fgrp$.  

Also we will be of course concerned with Kan extensions along groupoid maps, which we know exist whenever $\cat{C}$ satisfies the following 
\begin{defin}
A category $\cat{C}$ admits essentially finite (co)limits, if it admits (co)limits indexed by objects of $\fgrp$. 
\end{defin}

The simplest example of essentially finite groupoids is given by the classifying space $BG$ of a finite group $G$, while finite vector spaces are clearly closed under essentially finite (co)limits and, having duals, allow us to define the map $\ga{}$ of the previous section. It is in this context that one can easily check the failure of $\ga{}$ to satisfy functoriality. One would in fact expect that for composable maps $X\xto{f}Y\xto{g}Z$ of groupoids the diagram
\[
\xymatrix{
\ra{g}\ra{f}  \ar[d]_{\iso} \ar[r]^{\ra{g}\ga{f}} & \ra{g}\la{f} \ar[r]^{\ga{g}}  & \la{g}\la{f} \ar[d]^{\iso} \\
\ra{(gf)} \ar[rr]_{\ga{gf}}  && \la{(gf)} 
}
\]
commutes.

With reference to \cref{leftkanvec,rightkanvec,,nakavect}, we have instead:

\begin{ex} \label{counterexcn}
Let $C_{n}$ be the cyclic group of order $n$ with generator $g$, and consider the composition
\[
\{*\} \xto{s} BC_{n} \xto{t} \{*\}
\] 
Letting $k$ the trivial representation of the point in $\vect_{k}$, one gets a diagram
\[
\xymatrix{
Hom_{C_{n}}(k,Hom_{k}(k[C_{n}],k))  \ar[d]_{\iso} \ar[rr]^-{\ga{t}\circ\ra{t}\ga{s}} &&  k\otimes_{C_{n}} k[C_{n}]\tensor k   \\
Hom_{k}(k,k) \ar[rr]_-{\ga{ts}}  && k\tensor k \ar[u]_{\iso} 
}
\]
If $\phi\in Hom_{C_{n}}(k,Hom_{k}(k[C_{n}],k))$, so that $\psi=\phi(1)$ is constant on the elements of $C_{n}$ with value $\psi(g^{j})=\bar{\psi}$, then $\phi$ is taken by the horizontal morphism to
\[
n\cdot 1\tensor_{C_{n}} 1\tensor \bar{\psi}
\]
while via the lower path the factor $n$ does not arise, and one has just
\[
1\tensor_{C_{n}} 1\tensor \bar{\psi}
\]
It is clear then, that  such discrepancy can not be eliminated whenever the characteristic of $k$ divides $n$. In characteristic zero, instead, one can still hope to rescale the maps to avoid the problem.
\end{ex}
One should notice as well that the anomaly highlighted above varies locally at each connected component of the spaces involved, as explained in the following

\begin{ex} \label{counterexcncm}
Let $C_{m},C_{n}$ be the cyclic groups of order $m$ and $n$ with generators  $h$ and $g$, and consider the composition
\[
\{x,y\} \xto{s} BC_{m}\sqcup BC_{n} \xto{t} \{*\}
\] 
where $s$ sends $x,y$ respectively to the unique zero-cells of $BC_{m}$ and $BC_{n}$. We now get two maps
\[
\xymatrix{
Hom_{C_{m}}(k,Hom_{k}(k[C_{m}],k))\oplus Hom_{C_{n}}(k,Hom_{k}(k[C_{n}],k)) 
\ar@/^1pc/[d]^{\ga{t}\circ\ra{t}\ga{s}}  \ar@/_1pc/[d]_{\ga{ts}} \\
(k\otimes_{C_{m}} k[C_{m}] \tensor k ) \oplus ( k\otimes_{C_{n}} k[C_{n}] \tensor k )
%
}
\]
Given a pair of morphisms $(\phi_{m},\phi_{n})$ (sending $1$ to  functions $\psi_{m}(h^{i})=\bar{\psi}_{m}$, $\psi_{n}(g^{j})=\bar{\psi}_{n}$) we obtain via the rightmost map the value
\[
(m\cdot 1\tensor_{C_{m}} 1\tensor \bar{\psi}_{m},\enskip  n\cdot 1\tensor_{C_{n}} 1\tensor \bar{\psi}_{n})
\]
The other map, as before, ignores the contributions and gives
\[
(1\tensor_{C_{m}} 1\tensor \bar{\psi}_{m},\enskip 1\tensor_{C_{n}} 1\tensor \bar{\psi}_{n})
\]
As for the previous example, in non-zero characteristic there is no way to cure the problem.
\end{ex}

\subsection{Weights and the Nakayama map}

We now try make sense of the discrepancies appeared in \cref{counterexcn,counterexcncm} and define a new map $\na{f}$ out of $\ga{f}$, with similar properties.\\

First, in view of \cref{grpBC}, we can give the following definition of Kan extensions along maps of groupoids

\begin{defin}\label{kangrp}
If $f\colon X\to Y$ is a map between essentially finite groupoids, let
\[
\xymatrix{
P_{y}\ar[d]_{f_{y}} \ar[r]^{p} & X \ar[d]^{f} \\
\{*\} \ar[r]_{y} & Y \ultwocell<\omit>{\pi}
}
\]
be the homotopy fiber of $f$ at the point $y$.  
Given $V\colon X\to\cat{C}$, the value at $y$ of the left (resp. right) Kan extension $\la{f}V$ ($\ra{f}V$) along $f$, is  given  by the colimit (resp. limit) of $p^{*}V$.
\end{defin}

We can  define an endomorphism of the trivial representation $1$ on $P_{y}$ as follows

\begin{defin} \label{weightdelta}
Let $f\colon X\to Y$ be a map between essentially finite groupoids. On each connected component $\hat{P}\hookrightarrow P_{y}$ of the fiber
\[
\xymatrix{
P_{y}\ar[d]_{f_{y}} \ar[r]^{p} & X \ar[d]^{f} \\
\{*\} \ar[r]_{y} & Y \ultwocell<\omit>{\pi}
}
\]
at $y$, the (essentially) unique map $\{*\}\xto{\iota} \hat{P}$ defines an endomorphism 
\[
1\xto{\etar{\iota}} \ra{\iota}\iota^{*}1 \xto{\ga{\iota}} \la{\iota}\iota^{*}1 \xto{\epsl{\iota}} 1
\]
of the trivial representation $1$ on $\hat{P}$.  
The collection of such maps gives then an endomorphism of $1\in\cat{C}^{P_{y}}$
\[
\delta_{f} \colon1\to 1
\]
and therefore a natural endomorphism of $p^{*}$
\[
p^{*}V\iso p^{*}V\tensor1\xto{id\tensor\delta_{f}} p^{*}V\tensor 1\iso p^{*}V
\]
which we shall denote again by $\delta_{f}$.
\end{defin}

Recall that for an object $X$ in a symmetric monoidal category with duals $\cat{C}$, the dimension of $X$ is the endomorphism 
\[
1_{\cat{C}} \xto{co} X\tensor X^{d} \xto{ev}  1_{\cat{C}}
\]
One can easily prove that, at each connected component  $\hat{P}\hookrightarrow P_{y}$, $\delta_{f}$ is the multiplication by the dimension of $\la{\iota}1$.  
Placing us back in the situation of \cref{counterexcncm}, the endomorphism $\delta_{t}$ associated to  $t\colon BC_{m}\sqcup BC_{n}\to\{*\}$ multiplies the vectors relative to $BC_{j}$ ($j=m,n$) exactly by $j$. 

It makes sense, then, to  give the following

\begin{defin}[Nakayama map]\label{naka} 
Let $f\colon X\to Y$ be a map of essentially finite groupoids, and assume that $\delta_{f}$ is invertible. \\
Define a natural transformation, the \emph{Nakayama map},
\[
\na{f}\colon\ra{f}\to \la{f}
\]
 as
\[
\ra{f}V\xto{\ga{f}} \la{f}V  \xto{\la{f}\delta_{f}^{-1}} \la{f}V
\]
\end{defin}

Notice that $\na{f}$ could  as well be defined as
\[
\ra{f}V \xto{\ra{f}\delta_{f}^{-1}} \ra{f}V\xto{\ga{f}} \la{f}V 
\]

We can now show that  most  properties of  $\ga{f}$, proved in the previous section, are satisfied by $\na{f}$  as well.

\begin{prop}  \label{nakamate}
Let
\[
\xymatrix@!0{
& P \ar[ld]_{p} \ar[rd]^{q} \\
M \ar[rd]_{g}   &&  N  \ar[ld]^{h}   \\
& Y \uutwocell<\omit>{\pi}
}
\]
be a homotopy pullback of essentially finite groupoids. Let $\cat{C}$ be a symmetric monoidal category with duals closed under essentially finite (co)limits, and assume that $\delta_{p},\delta_{h}$ are invertible; then, the following diagram commutes
\[
\xymatrix{
\ra{p}q^{*}V  \ar[r]^{\na{p}}  &    \la{p}q^{*}V \ar[d]^{{\lma{(\pi^{-1})}}} \\
g^{*}\ra{h}V \ar[u]^{\rma{\pi}} \ar[r]_{g^{*}\na{h}} & g^{*}\la{h}V
}
\]
\proof
We can split the square as 
\[
\xymatrix{
\ra{p}q^{*}V  \ar[r]^{\ga{p}}  &    \la{p}q^{*}V \ar[d]^{{\lma{(\pi^{-1})}}} \ar[r]^{\la{p}\delta_{p}^{-1}}  &    \la{p}q^{*}V \ar[d]^{{\lma{(\pi^{-1})}}}  \\
g^{*}\ra{h}V \ar[u]^{\rma{\pi}} \ar[r]_{g^{*}\ga{h}} & g^{*}\la{h}V \ar[r]_{g^{*}\la{h}\delta_{h}^{-1}} & g^{*}\la{h}V
}
\]

Since we are dealing with a homotopy pullback, which is is exact by \cref{grpBC},  the square on the left commutes by \cref{gamate}.  
We have to show that the multipliers $\delta$ on the right do not affect commutativity. This is a consequence of the above definition of Kan extensions: in fact at each point both $\delta$'s arise in the same way, from the following homotopy pullback
\[
\xymatrix{
P' \ar[r] \ar[d] & P \ar[r] \ar[d] & N\ar[d] \\
\{*\} \ar[r] & M \ar[r] &  Y
}
\]
\endproof
\end{prop}

\begin{rem} \label{nakahomtopyrem}
Notice that, as a particular case, the above result states that Nakayama maps obtained from different choices of adjoints to a restriction functor are canonically conjugated.
\end{rem}

In analogy with \cref{pullpushga} we have

\begin{thm} \label{pullpushkangrp} 
Let
\[
\xymatrix@!0{
& P \ar[ld]_{p} \ar[rd]^{q} \\
M \ar[rd]_{g}   &&  N  \ar[ld]^{h}   \\
& Y \uutwocell<\omit>{\pi}
}
\]
be a homotopy pullback of essentially finite groupoids and $\cat{C}$ be a symmetric monoidal category with duals closed under essentially finite (co)limits, such that $\delta_{p},\delta_{h}$ are invertible. 
Let $\bar{\pi}\colon \la{g}\la{p}\to\la{h}\la{q}$ be the natural isomorphism induced by $\pi$. Then, the following diagram commutes 
\[
\resizebox{\textwidth}{!}{
\xymatrix@C+1pc{
\la{g}g^{*} \ar[rrd]_{\etar{h}} \ar[r]^{\la{g}\etar{p}} &  \la{g}\ra{p}p^{*}g^{*} \ar[r]^{\la{g}\ra{p}\pi}  & \la{g}\ra{p}q^{*}h^{*}  \ar[r]^{\la{g}\na{p}} &  \la{g}\la{p}q^{*}h^{*} \ar[r]^{\bar{\pi}} &  \la{h}\la{q}q^{*}h^{*} \ar[r]^{\la{h}\epsl{q}} &  \la{h}h^{*} \\
&& \ra{h}h^{*}\la{g}g^{*} \ar[r]_{\na{h}} &  \la{h}h^{*}\la{g}g^{*} \ar[rru]_{\la{h}h^{*}\epsl{g}} &&
}
}
\]
\proof
Being the square a homotopy pullback, we can take advantage of \cref{grpBC} and \cref{nakamate} above. The proof then proceeds exactly as in \cref{pullpushga}.
\endproof
\end{thm}

Finally, as $\ga{}$ in \cref{multiplyga}, also $\na{}$ can be seen to define a ``monoidal transformation'' between the ``monoidal functors'' given by left and right Kan extensions:

\begin{thm} \label{multiplynaka}
Let $f\colon M\to X$ and $g\colon N\to Y$ be maps,  
and $\cat{C}$ be a symmetric monoidal category with duals closed under essentially finite (co)limits, such that $\delta_{f},\delta_{g}$ are invertible. 
 Then for any ${V\in \cat{C}^{M}},{W\in\cat{C}^{N}}$ 
\[
\xymatrix{
\ra{(f\times g)}(V\etensor W) \ar[d]_{\na{f\times g}}  \ar[r]^{\ra{\mu}} &  \ra{f}V\etensor \ra{g}W \ar[d]^{\na{f}\etensor\na{g}} \\
\la{(f\times g)}(V\etensor W) \ar[r]_{\la{\mu}} & \la{f}V\etensor \la{g}W
}
\]
commutes.
\proof
By \cref{multiply}, $\la{\mu}$ and $\ra{\mu}$ are the left and right mates for  diagrams of the form 

\[
\xymatrix@C+2pc @R+1pc{
\cat{C}^{X}\times \cat{C}^{Y} \ar[d]_{f^{*}\times id} \ar[r]^{p_{X}^{*}\times p_{Y}^{*}}    & \cat{C}^{X\times Y} \times \cat{C}^{X\times Y} \ar[r]^{\tensor} \ar[d]|{(f\times id)^{*}\times id} & \cat{C}^{X\times Y} \ar[d]^{(f\times id)^{*}} \\
\cat{C}^{M} \times \cat{C}^{Y}  \ar[r]^{p_{M}^{*}\times p_{Y}^{*}}   &  \cat{C}^{M\times Y} \times  \cat{C}^{X\times Y} \ar[r]^-{\tensor (f\times id)^{*}}  &   \cat{C}^{M\times Y}  \ultwocell<\omit>{\mu} \\
}
\]

The left sided square commutes  with $\na{}$ by \cref{nakamate}. For the right sided one, we can apply the same ``splitting'' trick as in \cref{nakamate} and, by \cref{gamate}, concentrate only on the factors $\delta$; since the only factor appearing is $\delta_{(f\times id)}$, which agrees with tensor product, commutativity is verified.
\endproof
\end{thm}

\subsection{Nakayama Categories}

We can finally collect the results obtained so far, and characterise the kind of categories we should feed our quantization with. 
It is clear that the main issue was functoriality of the map $\ga{}$. In view of the linear case, we have introduced the weights $\delta$ and the new map $\na{}$. This in fact can be seen (\cref{nakavect}) to solve our problem over $\vect$, in addition to recovering and formalizing results and definitions from the previous literature \cite{bd2,fhlt,mo1,mo2}. 
We will therefore ask our target category $\cat{C}$ to guarantee existence and functoriality of the Nakayama map  as follows:

\begin{defin} \label{nakacontext}
We say that a category $\cat{C}$ is a \emph{Nakayama Category} if
\begin{itemize}
\item[i)] $\cat{C}$ has essentially finite colimits.
\item[ii)] $\cat{C}$ is symmetric monoidal.
\item[iii)] $\cat{C}$ has duals.
\item[iv)] For any map $X\xto{f}Y$ in $\fgrp$, the endomorphism $\delta_{f}$ is invertible.
\item[v)] For any composable pair of maps $X\xto{f}Y\xto{g}Z$ in $\fgrp$ the diagram
\[
\xymatrix@!=5pt{
& \la{g}\ra{f} \ar[rd]^{\la{g}\na{f}} & \\
\ra{g}\ra{f} \ar[rd]_{\ra{g}\na{f}} \ar[ur]^{\na{g}} \ar[rr]^{\na{gf}} && \la{g}\la{f} \\
& \ra{g}\la{f} \ar[ru]_{\na{g}} &
}
\]
commutes.
\end{itemize}
\end{defin}

Let us further motivate the above definition. Essentially, we would like to apply the machinery developed in the previous sections to restriction functors $f^{*}\colon \cat{C}^{Y}\to\cat{C}^{X}$ and their adjoints, given by Kan extensions.  
Since we are dealing with groupoids points $ii)$ and $iii)$ ensure, by \cref{fundualcat}, that the categories $\cat{C}^{X}$ we will consider are symmetric monoidal with duals.  
The restriction to essentially finite groupoids in point $i)$ is motivated by concrete examples of categories with duals, such as vector spaces, where infinite colimits are not dualizable objects. 
Point $i)$ ensures the existence of left Kan extensions (along maps of essentially finite groupoids) and therefore of adjoint pairs of functors $\la{f}\lad f^{*}$ between $\cat{C}^{M}$ and $\cat{C}^{X}$, for any ${M\xto{f}X}$. 
Being the restriction functors (symmetric) monoidal we can now apply \cref{rightiffleft}, and conclude that also right Kan extensions exist. Points $i)-v)$, altogether, will allow us to use to results of the previous pages so to define the quantization functors on objects and, via the Nakayama map, on morphisms.  
In particular, point $v)$ will be necessary in \cref{thmsum,thmprod}, to
verify functoriality of $\sam$ and $\pro$ by taking advantage of \cref{pullpushkangrp}. 

\begin{rem} \label{nakacontextrem2}
In view of the above comments and \cref{rightiffleft}, one can reformulate \cref{nakacontext}, by asking for (essentially finite) limits instead of colimits. The two formulations are of course equivalent.
\end{rem}

\begin{rem}
Notice that condition $v)$ of \cref{nakacontext}, and in general the definition of $\na{f}$, does not depend on the choice of representatives for the Kan extensions along $f$. In fact, by \cref{nakahomtopyrem}, Nakayama maps relative to isomorphic pairs of functors are conjugated.
\end{rem}

\section{Quantization}\label{secquant}

The techniques developed in the previous sections will  now be used to define, in terms of Kan extensions, two monoidal functors
\begin{align*}
\sam \colon & \fami{\cat{C}} \to \cat{C}  \\
\pro \colon & \fami{\cat{C}} \to \cat{C} 
\end{align*}
from a category of representations of essentially finite groupoids valued in any Nakayama category $\cat{C}$.

We'll see in \cref{secsumisopro} that the two functors are canonically isomorphic, thus implying invertibility of the Nakayama morphism $\na{}$.

\subsection{A monoidal category of local systems} \label{secfam}

Here below we will briefly describe the source category $\fami{\cat{C}}$ of the quantization functors.  

In the general case, this is an $(\infty,n)$-category $\famic{n}{\cat{C}}$ ``of spans of $\infty$-groupoids''. A typical object is an $\infty$-groupoid $X$ together with a functor $X\to\cat{C}$ to the target $(\infty,n)$-category $\cat{C}$. The $1$-morphisms are given by spans ${X\leftarrow M\to Y}$ along with a filling cell making the resulting diagram commute. A $2$-morphism from $M$ to $M'$ is given by a span of $\infty$-groupoids $M\leftarrow N\to M'$, with a filling two-cell; $k$-morphisms until level $n$ are similarly given by spans between spans, while higher ($k>n$) $k$-morphisms are given by (higher) homotopies between maps of $\infty$-groupoids, compatible with the underlying  diagrams.  
A precise description of such $(\infty,n)$-category $\famic{n}{\cat{C}}$ has been given in \cite{ha2}; there the author uses the notation $\mathit{Span}_{n}(\cat{S};\cat{C})$, while we preferred to keep the one used in \cite{fhlt} and \cite{lu3}.  

In our situation, being $\cat{C}$ merely a category, we need to restrict our attention to (finite) homotopy $1$-types (i.e. essentially finite groupoids) in order to use the results of \cref{secnaka}. Also, we can avoid working with higher categories and just consider the homotopy category of $\fami{\cat{C}}$ as a source. By an abuse of notation, this will again be denoted by $\fami{\cat{C}}$, and we will construct functors
\[
\fami{\cat{C}}\to\cat{C}
\]
One can anyway check that the two functors lift to $\infty$-functors (where now the target will be the classifying diagram \cite{re} of $\cat{C}$).\\

An object of $\fami{\cat{C}}$ will be given, as already mentioned, by an essentially finite groupoid $X$, and a functor
\[
X\xto{V}\cat{C}
\]
Morphisms will be (equivalence classes of) diagrams
\[
\xymatrix@!0{
&  M \ar[ld]_{f} \ar[rd]^{g} \\
X \ar[rd]_{V}  &&   Y  \ar[ld]^{W}   \\
&  \cat{C}  \uutwocell<\omit>{\alpha}
}
\]
with composition  given by (an equivalence class of) homotopy pullback along the middle maps
\[
\xymatrix@!0{
&&  P \ar[ld]_{p} \ar[rd]^{q} \\
&M \ar[ld]_{f} \ar[rd]_{g}  &&   N  \ar[ld]^{h} \ar[rd]^{k}  \\
X \ar[rrdd]_{U} &&  Y \ar[dd]^{V} \uutwocell<\omit>{\pi} && Z  \ar[lldd]^{W} \\
&{} \utwocell<\omit>{\alpha} & & {} \utwocell<\omit>{\beta} \\
&& \cat{C}
}
\]
and the identity morphism is clearly given by
\[
\xymatrix@!0{
&  X \ar[ld]_{id} \ar[rd]^{id} \\
X \ar[rd]_{V}  &&   X  \ar[ld]^{V}   \\
&  \cat{C}  \uutwocell<\omit>{id}
}
\]

Cartesian product and external tensor product (\cref{exttensor}) endow $\fami{\cat{C}}$ with the structure of a symmetric monoidal category with duals.

Given objects $X\xto{V}\cat{C},Y\xto{W}\cat{C}$, their tensor product is
\[
X\times Y \xto{V\etensor W} \cat{C}
\]
The unit is provided by the unit of $\cat{C}$, viewed as a functor
\[
*\xto{1}\cat{C}
\]
from the one-point space.

Finally, the dual to $X\xto{V}\cat{C}$ is given by $X\xto{V^{d}}\cat{C}$ as from \cref{fundualcat}, with coevaluation and evaluation morphism

\[
\xymatrix@!0{
&  X \ar[ld] \ar[rd]^{\Delta} \\
{*} \ar[rd]_{1}  &&   X\times X  \ar[ld]^{V\etensor V^{d}}   \\
&  \cat{C}  \uutwocell<\omit>{co}
}
\qquad
\xymatrix@!0{
&  X \ar[ld]_{\Delta} \ar[rd] \\
X\times X \ar[rd]_{V^{d}\etensor V}  &&   {*}  \ar[ld]^{1}   \\
&  \cat{C}  \uutwocell<\omit>{ev}
}
\]

\subsection{The quantization functors}

It  is not difficult, at this point, to define the functors $\sam$ and $\pro$. As already mentioned in the introduction the value on objects will be given by (co)limits, while the map $\na{}$ will be used to define the functors on morphisms. The results of \cref{nakakantensor} and the compatibility with $\na{}$ proved in \cref{multiplynaka} will make them into symmetric monoidal functors.

In the following we will work up to the obvious equivalence $\cat{C}\simeq\cat{C}^{*}$. 
For a space $X$, denote by $x$ the unique map $X\to*$ to the point.\\

\begin{thm} \label{thmsum}\index{sum@$\sam$}
Let $\cat{C}$ be a Nakayama category. Then, there exists a symmetric monoidal functor 
\[
\sam\colon \fami{\cat{C}}\to\cat{C}
\]
\proof
Given an object $V\colon X\to\cat{C}$ of $\fami{\cat{C}}$, we let
\[
\sam V= \la{x}V
\]
be the left Kan extension of $V$ along $x$, i.e. its colimit.

To a morphism
\[
\xymatrix@!0{
&  M \ar[ld]_{f} \ar[rd]^{g} \\
X \ar[rd]_{V}  &&   Y  \ar[ld]^{W}   \\
&  \cat{C}  \uutwocell<\omit>{\alpha}
}
\]
$\sam$ associates the map $\sam(\alpha)$
\[
\resizebox{\textwidth}{!}{
$
\la{x}V\xto{\la{x}\etar{f}} \la{x}\ra{f}f^{*}V\xto{\la{x}\ra{f}\alpha}\la{x}\ra{f}g^{*}W\xto{\la{x}\na{f}}\la{x}\la{f}g^{*}W\xto{\iso}\la{y}\la{g}g^{*}W\xto{\la{y}\epsl{g}}\la{y}W
$}
\]
It is clear that on identity morphisms, i.e. squares
\[
\xymatrix@!0{
&  X \ar[ld]_{id} \ar[rd]^{id} \\
X \ar[rd]_{V}  &&   X  \ar[ld]^{V}   \\
&  \cat{C}  \uutwocell<\omit>{\alpha}
}
\]
the above assignment reduces to the identity map.

Now let's consider a composition of morphisms in $\fami{\cat{C}}$, that is a diagram
\[
\xymatrix@!0{
&&  P \ar[ld]_{p} \ar[rd]^{q} \\
&M \ar[ld]_{f} \ar[rd]_{g}  &&   N  \ar[ld]^{h} \ar[rd]^{k}  \\
X \ar[rrdd]_{U} &&  Y \ar[dd]^{V} \uutwocell<\omit>{\pi} && Z  \ar[lldd]^{W} \\
&{} \utwocell<\omit>{\alpha} & & {} \utwocell<\omit>{\beta} \\
&& \cat{C}
}
\]
where the upper square is (by definition of $\fami{\cat{C}}$) a homotopy pullback of $1$-types. 
The outer square and the two diagrams labeled by $\alpha$ and $\beta$ give rise to $1$-morphisms $\sam( \beta \pi \alpha ),\sam(\alpha),\sam( \beta )$ of $\cat{C}$. I claim that the first is the composition of the other two.  

Noting that the isomorphism $\la{x}\la{f}\la{p}\to\la{z}\la{k}\la{q}$ is the composition 
\[
\la{x}\la{f}\la{p}\xto{\iso}\la{y}\la{g}\la{p}\xto{\la{y}\bar{\pi}}\la{y}\la{h}\la{q}\xto{\iso}\la{z}\la{k}\la{q}
\]
a simple analysis of the resulting diagram using the functoriality hypothesis and \cref{pullpushkangrp} proves the claim.  

It remains to show that $\sam$ is actually a symmetric monoidal functor. 
The isomorphism $\la{\mu}$ of \cref{multiply} provides an isomorphism
\[
\la{\mu}\colon\sam(V\etensor V')\to\sam(V)\tensor\sam(V')
\]
Naturality of $\la{\mu}$, i.e. commutativity of
\[
\xymatrix{
\sam(V\etensor V')\ar[d]_{\sam(\alpha\etensor\alpha)}\ar[r]^-{\la{\mu}} & \sam(V)\tensor\sam(V') \ar[d]^{\sam(\alpha)\tensor \sam(\alpha')} \\
\sam(W\etensor W')\ar[r]_-{\la{\mu}} & \sam(W)\tensor\sam(W') 
}
\]
for morphisms $\alpha,\alpha'$ and the corresponding tensor product
\[
\xymatrix@!0{
&  M \ar[ld]_{f} \ar[rd]^{g} \\
X \ar[rd]_{V}  &&   Y  \ar[ld]^{W}   \\
&  \cat{C}  \uutwocell<\omit>{\alpha}
}
\qquad
\xymatrix@!0{
&  M' \ar[ld]_{f'} \ar[rd]^{g'} \\
X' \ar[rd]_{V'}  &&   Y'  \ar[ld]^{W'}   \\
&  \cat{C}  \uutwocell<\omit>{\alpha'}
}
\qquad
\xymatrix@!=7pt{
&  M \times M' \ar[ld]_{f\times f'} \ar[rd]^{g\times g'} \\
X\times X' \ar[rd]_{V\etensor V'}  &&   Y\times Y'  \ar[ld]^{W\etensor W'}   \\
&  \cat{C}  \uutwocell<\omit>{\alpha\etensor\alpha'}
}
\]
follows from \cref{multiplycomp,multiplyhomo,multiplyetaeps,multiplynaka}.  

The unit isomorphism 
\[
\sam(1)\to1
\]
is just the identity, being $\sam(1)$ simply the Kan extension along the identity map $*\to*$.

Finally  associativity, unitality and commutativity of $\sam$ follow by \cref{multiplyassociativity,multiplycommute}.
\endproof
\end{thm}

As the reader might have noticed, the arguments used in the previous theorem can all be dualized. It is natural then to consider right Kan extensions in place of left ones, and try to define a ``limit'' version of $\sam$. To further stress this duality, recall from \cref{nakacontextrem2} that we can define a Nakayama category as having (essentially finite) limits instead of colimits.

A similar proof  gives

\begin{thm} \label{thmprod}
Let $\cat{C}$ be a Nakayama category. 
Then, there exists a symmetric monoidal functor
\[
\pro\colon \fami{\cat{C}}\to\cat{C}
\]
\proof
The proof is dual to the one of \cref{thmsum}. 
This time, $\pro$ associates to $V\colon X\to\cat{C}$ the right Kan extension $\ra{x}V$ (i.e. the limit of $V$).  

A morphism 
\[
\xymatrix@!0{
&  M \ar[ld]_{f} \ar[rd]^{g} \\
X \ar[rd]_{V}  &&   Y  \ar[ld]^{W}   \\
&  \cat{C}  \uutwocell<\omit>{\alpha}
}
\]
is taken to the map
\[
\resizebox{\textwidth}{!}{
$
\ra{x}V\xto{\ra{x}\etar{f}} \ra{x}\ra{f}f^{*}V \xto{\iso}\ra{y}\ra{g}f^{*}V  \xto{\ra{y}\ra{g}\alpha}\ra{y}\ra{g}g^{*}W  \xto{\ra{y}\na{g}}\ra{y}\la{g}g^{*}W  \xto{\ra{y}\epsl{g}}\la{y}W
$}
\]
The multiplication isomorphism
\[
\pro(V\etensor V') \to \pro(V)\tensor \pro(V')
\]
is clearly given by the isomorphism $\ra{\mu}$ of \cref{multiply}.
\endproof
\end{thm}

\subsection{Siamese twins}\label{secprod} \label{secsumisopro}

We can finally conclude by comparing the two functors just defined. Recalling from \cref{multiplynaka} the ``monoidal'' behaviour of the transformation $\na{}$, the following comes with no surprise

\begin{prop}\label{nakamontran}
Let $\cat{C}$ be a Nakayama category. Then the Nakayama map defines a monoidal transformation
\[
\na{}\colon\pro\to\sam
\]
\proof
First, we check that the maps
\[
\na{x}\colon \pro(V)=\ra{x}V\to\la{x}V=\sam(V)
\]
define a natural transformation. For any span
\[
\xymatrix@!0{
&  M \ar[ld]_{f} \ar[rd]^{g} \\
X \ar[rd]_{V}  &&   Y  \ar[ld]^{W}   \\
&  \cat{C}  \uutwocell<\omit>{\alpha}
}
\]
we obtain in fact by the functoriality assumption and \cref{nakamate}, a commutative diagram
\[
\resizebox{\textwidth}{!}{
\xymatrix{
\ra{x}V \ar[r]^{\ra{x}\etar{f}} \ar[dd]_{\na{x}} & \ra{x}\ra{f}f^{*}V \ar[r]^{\ra{x}\ra{f}\alpha} &  \ra{x}\ra{f}g^{*}W   \ar[dd]_{\na{x}}   \ar[ddr]_{\na{xf}} \ar[r]^{\rma{id}} &  \ra{y}\ra{g}g^{*}W  \ar[ddr]^{\na{yg}}   \ar[r]^{\ra{y}\na{g}} \ar@{}[dd]|{\cref{nakamate}} & \ra{y}\la{g}g^{*}W \ar[dd]^{\na{y}} \ar[r]^{\ra{y}\epsl{g}} &  \ra{y}W \ar[dd]^{\na{y}} \\
&&&    &&  \\
\la{x}V \ar[r]_{\la{x}\etar{f}} & \la{x}\ra{f}f^{*}V \ar[r]_{\la{x}\ra{f}\alpha} &  \la{x}\ra{f}g^{*}W  \ar[r]_{\la{x}\na{f}}   &  \la{x}\la{f}g^{*}W      \ar[r]_{\lma{id}} & \la{y}\la{g}g^{*}W  \ar[r]_{\la{y}\epsl{g}} &  \la{y}W
}
}
\]
exhibiting $\na{}$ as a natural transformation. 

The fact that $\na{}$ is monoidal, follows by simply applying \cref{multiplynaka} to
\[
\xymatrix{
\ra{(x\times x')}(V\etensor V') \ar[d]^{\ra{\mu}} \ar[r]^{\na{}} &  \la{(x\times x')}(V\etensor V') \ar[d]^{\la{\mu}} \\
\ra{x}V\tensor \ra{x'}V' \ar[r]_{\na{}\tensor\na{}} & \la{x}V\tensor \la{x'}V'
}
\]
and trivially from
\[
\xymatrix{
& 1 \ar[dl]_{id} \ar[dr]^{id} \\
\ra{id}1 \ar[rr]_{\na{}} && \la{id}1
}
\]
So that $\na{}$ is compatible with the monoidal structure of $\pro$ and $\sam$.
\endproof
\end{prop}

We can  summarize the previous section and the above result in a single theorem:

\begin{thm}\label{sumisopro}
Let $\cat{C}$ be a Nakayama category, then there exist canonically isomorphic symmetric monoidal functors
\begin{align*}
\pro \colon & \fami{\cat{C}} \to \cat{C} \\
\sam \colon & \fami{\cat{C}} \to \cat{C} 
\end{align*}
respectively defined by right and left Kan extensions of local systems.
\proof
We already know that the Nakayama map defines a monoidal transformation ${\pro\to\sam}$. Being $\fami{\cat{C}}$ a category with duals, it follows from \cref{montraiso} that $\na{}$ is actually an isomorphism.
\endproof
\end{thm}

As a consequence, we obtain that in any Nakayama category $\cat{C}$ limits and colimits are canonically isomorphic, via either $\na{}$ or $\ga{}$, so that the assumption in \cite{fhlt,mo1,mo2} is now a corollary. Of course, this extends to an isomorphism between left and right Kan extensions, generalising a well-known fact from representation theory.

\begin{cor} \label{nakaiso}
Let $\cat{C}$ be a Nakayama category. Then, for any map ${X\xto{f}Y}$ of essentially finite groupoids, the Nakayama map defines a canonical isomorphism
\[
\ra{f}V\to\la{f}V
\]
for any $V\colon X\to \cat{C}$.
\proof
We know that the statement is true when $f$ is the terminal map to the point, so that the right and left Kan extensions are actually the limit and colimit functors. From \cref{kangrp} it follows that for any map $X\xto{f}Y$, the value of $\na{f}$ at any point $y\in Y$ is an isomorphism as well.
\endproof
\end{cor}

\section{Linear representations of groupoids}\label{secvect}

We conclude by rapidly considering the results of the previous sections in the context of finite vector spaces.  
While in this situation the functors $\sam$ and $\pro$ themselves do not give rise to really interesting theories, the techniques used to build them and in particular \cref{nakaiso} specialize to relevant known results, which appear now under a different light.\\

Recall that a representation $V$ of a group $G$ is the same thing as a functor $V\colon BG\to\vect$ from the classifying space $BG$ of $G$. Similarly, a $\vect$-valued representation of a groupoid $X$ is a functor $V\colon X\to\vect$. 

Also notice that any essentially finite groupoid $X$ is equivalent to a finite disjoint union $\coprod BG_{i}$, with the $G_{i}$'s finite groups. Therefore we will simplify notations (and calculations) by assuming that  $\map{x}{y}$ is empty if $x\neq y$, and denote $A_{x}=\map{x}{x}$.

In the following  $\vect$ will mean finite (hence dualizable) vector spaces over a field $k$ of characteristic zero.

The usual presentation of induction and coinduction of group representations in terms of tensor and hom spaces then naturally generalises to essentially finite groupoids:

\begin{thm}  \label{leftkanvec}
Let $f\colon X\to Y$ be a morphism between essentially finite groupoids, and let $V\colon X\to \vect$ be a representation of $X$. Then
\begin{itemize}
\item
The value of $\la{f}V$ at $y\in Y$ is
\[
\la{f}V(y)= \bigoplus_{x|fx=y}  k[A_{y}]\tensor_{k[A_{x}]} V(x)
\]
\item
At $x\in X$, the unit $\eta\colon V(x)\to f^{*}\la{f}V(x)=\la{f}V(fx)$ is the map
\[
v\mapsto 1\tensor_{k[A_{x}]} v \in k[A_{y}]\tensor_{k[A_{x}]} V(x)
\]
\item The counit $\eps\colon\la{f}f^{*}W(y) \to W(y)$ is the map induced by the assignments
\[
k[A_{y}]\tensor_{k[A_{x}]} W(y) \ni g\tensor_{k[A_{x}]} w \mapsto gw 
\]
where $g\in A_{y}$
\end{itemize}
\end{thm}

\begin{thm} \label{rightkanvec}
Let $f\colon X\to Y$ be a morphism between essentially finite groupoids, and let $V\colon X\to \vect$ be a representation of $X$. Then
\begin{itemize}
\item
The value of $\ra{f}V$ at $y\in Y$ is
\[
\ra{f}V(y) =\bigoplus_{x|fx=y}  Hom_{k[A_{x}]}(k[A_{y}],V(x))
\]
\item
At $y\in Y$, the unit $\eta\colon W(y)\to \ra{f} f^{*}W(y)$ is the map
\[
w\mapsto \sum_{x|fx=y} \phi_{w,x}
\]
where $\phi_{w,x}(g)= gw$ for $g\in A_{y}$.
\item 
The counit $\eps\colon f^{*}\ra{f}V(x)=\ra{f}V(fx) \to V(x)$ is the map
\[
\sum_{x'|fx'=fx} \phi_{x'} \mapsto \phi_{x}(1)
\]
\end{itemize}
\end{thm}

One can now compute the Nakayama morphism for a map $f$ of essentially finite groupoids, according to the above presentations of Kan extensions.

\begin{thm}\label{nakavect} \label{gammavect}
Let $f\colon X\to Y$ be a morphism between essentially finite groupoids, and $V\colon X\to \vect$  a representation of $X$. 
For $x\in X$ let $f(A_{x})$ be the image of $A_{x}$ in $A_{fx}$, and denote by $K_{x}$ the kernel of the induced map $A_{x}\xto{f}A_{fx}$. 
Then
\begin{itemize}
\item The map $\ga{f}$ at $y\in Y$
\[
\ga{f,y} \colon \bigoplus_{x|fx=y}  Hom_{k[A_{x}]}(k[A_{y}],V(x)) \to  \bigoplus_{x|fx=y}  k[A_{y}]\tensor_{k[A_{x}]} V(x)
\]
is induced by the assignments
\[
Hom_{k[A_{x}]}(k[A_{y}],V(x)) \ni \phi_{x} \mapsto \frac{1}{|f(A_{x})|} \sum_{g\in A_{y}} g^{-1}\tensor\phi_{x}(g)
\]
\item 
The map $\delta_{f}\colon \la{f}V(y)\to\la{f}V(y)$ at $y$ is given at each summand by multiplication by 
\[
|K_{x}|
\]
\item The Nakayama map at $y$ is induced at each summand by the map
\[
\phi_{x} \mapsto \frac{1}{|A_{x}|} \sum_{g\in A_{y}} g^{-1}\tensor\phi_{x}(g)
\]
\end{itemize}
\end{thm}

As we already anticipated, we have that

\begin{thm}\label{vectnakacont}
The category of finite dimensional vector spaces over a zero-characteristic field is a Nakayama category.
\proof
We know that the category of finite dimensional vector spaces is symmetric monoidal with duals, and it has essentially finite colimits; if the characteristic of the field is zero, then $\delta_{f}$ is invertible for any $f$. The functoriality of the Nakayama map follows from a tedious but not difficult computation.
\endproof
\end{thm}

With \cref{gammavect,vectnakacont} at hand, one can easily check that in the case of trivial representations of groupoids the above recovers the values in \cite{mo1}. 
In particular, invertibility of the Nakayama morphism and ambidexterity of the adjunction given by restriction appear now as a consequence rather than an ingredient to quantization. \cref{nakaiso} in fact specializes to

\begin{thm} \label{frobpair}
Let $k$ be a zero-characteristic field, $H\xto{f}G$ a morphism of finite groups, and $coInd_{H}^{G},Ind_{H}^{G}\colon Rep_{k}(H)\to Rep_{k}(G)$ the coinduction and induction functors. Then
\begin{itemize}
\item[i)] There is a canonical natural isomorphism $\nu\colon coInd_{H}^{G}\to Ind_{H}^{G}$, the Nakayama isomorphism
\item[ii)] Induction is both a left and a right adjoint to restriction
\end{itemize} 
\end{thm}


\bigskip








\begin{thebibliography}{99}

\bibitem[At]{at} M. Atiyah, \emph{Topological quantum field theories}, Publ. Math. IHES 68, 1988, 175-186

\bibitem[BD1]{bd} J. Baez, J. Dolan, \emph{Higher-dimensional algebra and topological quantum field theory}, J. Math. Phys. 36 (11), 1995, 6073-6105

\bibitem[BD2]{bd2} J. Baez, J. Dolan, \emph{From finite sets to Feynman diagrams}, Mathematics Unlimited-2001 And Beyond. Springer, New York, 2001, 29-50

\bibitem[Ben]{ben} D. J. Benson, \emph{Representations and cohomology, Vol. 1: Basic representation theory of finite groups and associative algebras}, Cambridge Studies in Advanced Mathematics 30, Cambridge University Press, 1991

\bibitem[CS]{cs} D. Calaque, C. Scheimbauer, \emph{A note on the $(\infty,n)$-category of cobordisms}, 2015, preprint  \href{http://arxiv.org/abs/1509.08906}{arXiv:1509.08906 [math.AT]}

\bibitem[CD]{cd} D.-C. Cisinski, F. D\'eglise, \emph{Triangulated categories of mixed motives}, 2012, preprint \href{http://arxiv.org/abs/0912.2110}{arXiv:0912.2110 [math.AG]}

\bibitem[FHLT]{fhlt} D.S. Freed, M.J. Hopkins, J. Lurie, C. Teleman, \emph{Topological Quantum Field Theories from Compact Lie Groups}, A celebration of the mathematical legacy of Raoul Bott, CRM Proc. Lecture Notes, vol. 50, AMS, 2010, 367-403

\bibitem[FHM]{fhm} H. Fausk, P. Hu, J. P. May, \emph{Isomorphisms between left and right adjoints}, Theory and Applications of Categories 11, 2003, 107-131

\bibitem[Ha1]{ha1} R. Haugseng, \emph{The Becker-Gottlieb transfer is functorial}, 2013, preprint \href{http://arxiv.org/abs/1310.6321}{arXiv:1310.6321 [math.AT]}

\bibitem[Ha2]{ha2} R. Haugseng, \emph{Iterated spans and classical topological field theories},  2014, preprint \href{http://arxiv.org/abs/1409.0837}{arXiv:1409.0837 [math.AT]}


\bibitem[HeL]{hel} G. Heuts, J. Lurie, \emph{Ambidexterity}, in Topology and Field Theories, Contemporary Mathematics, Volume 613, AMS, 2014

\bibitem[HoL]{hol} M. Hopkins, J. Lurie, \emph{Ambidexterity in $K(n)$-Local Stable Homotopy Theory}, 2013, available \href{http://www.math.harvard.edu/~lurie/papers/Ambidexterity.pdf}{online}

\bibitem[JT]{jt} A. Joyal, M. Tierney, \emph{Notes on simplicial homotopy theory}, CRM Barcelona,  2008

\bibitem[Lu1]{lu3} J. Lurie, \emph{On the Classification of Topological Field Theories}, Current Developments in Mathematics, Volume 2008 (2009), 129-280

\bibitem[Lu2]{lu4} J. Lurie, \emph{Higher Algebra}, 2014, available \href{http://math.harvard.edu/~lurie/papers/higheralgebra.pdf}{online}

\bibitem[ML]{ml} S. Mac Lane, \emph{Categories for the Working Mathematician}, Graduate Texts in Mathematics 5, Springer, 1971

\bibitem[Mo1]{mo1} J.C. Morton, \emph{Two-Vector Spaces and Groupoids}, Applied Categorical Structures, vol 19, no. 4, 2011

\bibitem[Mo2]{mo2} J.C. Morton, \emph{Cohomological Twisting of $2$-Linearization and Extended TQFT}, J. Homotopy Relat. Struct. (2015) 10, 127-187

\bibitem[Re]{re} C. Rezk, \emph{A model for the homotopy theory of homotopy theory}, Trans. Amer. Math. Soc., 353(3), 2001, 973-1007

\end{thebibliography}
\end{document}